\documentclass{ijclclp}

\title{The planar Turán number of double star $S_{3,5}$}
\author{Dan-Dan Liu, Shou-Jun Xu* \\
        School of Mathematics and Statistics, Gansu Center for Applied Mathematics, Lanzhou University, Lanzhou, Gansu 730000, China}
\email{\texttt{liudd2023@lzu.edu.cn, shjxu@lzu.edu.cn}}
\usepackage{amsthm}

\usepackage{tipx}  
\DeclareUnicodeCharacter{030B}{\doubleacute}  
\usepackage{doi}
\usepackage[numbers]{natbib}

\begin{document}
\maketitle

\thispagestyle{firstpage}
\begin{abstract}
Given a graph $H$ and a positive integer $n$, the planar Turán number of $H$, denoted by $ex_\mathcal{P}(n, H)$, is the maximum number of edges in an $n$-vertex $H$-free planar graph. D.Ghosh, et al.initiated the topic of double stars $S_{k,l}$. Recently, Xu et al.[AIMS Mathematics, 2025, 10(1): 1628–1644] mentioned that $ex_\mathcal{P}(n, S_{3,5})$ is still unknown.  In this paper, we first establish that the planar Turán number $S_{3,5}$ satisfies $ex_\mathcal{P}(n, S_{3,5})$ $\leq$ $\tfrac{23n}{8}-\tfrac{9}{2}$ for all $n \geq 2$. The upper bound is tight when $n=12$.

\textbf{Keywords:} 
Planar Turán number, Double star, Extremal planar graph
\end{abstract}

\section{Introduction }
This paper restricts attention to finite simple graphs, where edges uniquely connect distinct vertex pairs.
 
 Let $G=(V(G), E(G))$ be a graph, where $V(G)$ and $E(G)$ are the vertex set and edge set of the graph $G$, respectively. A planar graph is a graph that can be drawn on a plane such that its edges are non-intersecting. We use
$N_G(v)$ to denote the set of vertices of $G$ adjacent to $v$, abbreviated without ambiguity as $N(v)$. The degree of a vertex $v$ is the number of vertices of $N(v)$, denoted by $d_G(x)$, abbreviated without ambiguity as $d(x)$ i.e. $|N(v)|$. Let $N_G[v]=N(v) \cup \{v\} $, abbreviated without ambiguity as $N[v]$. Let $v(G)$, $e(G)$, $\delta(G)$and $\Delta(G)$ denote the number of
vertices, the number of edges, the minimum degree and the maximum degree of $G$, respectively. For any
subset $S\subset V(G)$, the subgraph induced on $S$ is denoted by $G[S]$. We denote by $G\backslash S$ the
subgraph induced on $V(G)\backslash S$. If $S={v}$, we simply write $G\backslash v$. We use $E[S,T]$ to denote the edge set between S and T, let $e[S,T]$=$|E[S,T]|$, where $S,T$ are subsets of $V(G)$.

Let $H$ be a graph, and a graph is called $H$-free if it does not contain $H$ as a subgraph. A well-known problem in extremal graph theory is the Turán problem: What is the maximum number of edges in an $n$-vertex graph that does not contain $H$ as a subgraph? How can we construct such a graph? Is this graph unique? This problem defines the Turán number.

Given a graph $H$, the Turán number of $H$, denoted
by $ex(n, H)$, is the maximum number of edges in an $H$-free graph on $n$ vertices. In this paper, we focus on the planar Turán number. The planar Turán number of a graph $H$, denoted
by $ex_\mathcal{P}(n, H)$, is the maximum number of edges in an $H$-free planar graph on $n$ vertices.

Significant progress has been made in the research on the planar Turán number. Dowden \cite{Dowden2015ExtremalCP} in 2016 initiated the study of planar Turán-type problems. Dowden studied the planar Turán number of $C_4$ and $C_5$, where $C_k$ is a
cycle with $k$ vertices. Then Ghosh, Győri, Martin, Paulos and Xiao \cite{doi:10.1137/21M140657X} gave the exact value for $C_6$.
After that, Shi, Walsh and Yu \cite{shi2023planarturannumber7cycle}, Győri, Li and Zhou \cite{2023planarturannumbersevencycle}  gave the exact value for $C_7$. But the planar Turán number of $C_k$ remins undermined for $k\ge 8$. And Lan, Shi and Song \cite{lan2018extremalhfreeplanargraphs} gave a sufficient
condition for graphs with planar Turán number of $3n-6$. 
We refer the interested readers to
more results on paths, theta graphs and other graphs \cite{ Du2021PlanarTN,Fang2023ExtremalSR,Lan2019PlanarTN,lan2019extremalthetafreeplanargraphs}.

 The double star is an important structure in graph theory, but the special structure of the double star has not been considered on planar graphs before. Then in 2022, Győri, Martin, Paulos and Xiao \cite{2021arXiv211010515G} studied the topic for double stars as the forbidden graph. A $(k,l)$-star, denoted by $S_{k,l}$, is the graph obtained from an edge $uv$, and joining end vertices with $k$ and $l$ vertices, respectively. A $k$-$l$ edge refers to an edge whose endpoints have degrees $k$ and $l$, respectively. Győri, Martin, Paulos and Xiao gave the exact value for $ex_\mathcal{P}(n, S_{2,2})$ and $ex_\mathcal{P}(n, S_{2,3})$, and gave the upper bounds of $ex_\mathcal{P}(n, S_{2,4})$, $ex_\mathcal{P}(n, S_{2,5})$, $ex_\mathcal{P}(n, S_{3,3})$, $ex_\mathcal{P}(n, S_{3,4})$. In 2024, Xu \cite{XU2024326} $et$ $al$. improved the upper bound of $ex_\mathcal{P}(n, S_{2,5})$. In 2025, Xu \cite{XU2025114571,XU11} $et$ $al$. improved the upper bounds of $ex_\mathcal{P}(n, S_{2,4})$ and $ex_\mathcal{P}(n, S_{3,4})$. Meanwhile, Xu $et$ $al$. mentioned that $ex_\mathcal{P}(n, S_{3,5})$ is still unknown. Then in this paper, we establish an upper bound for $ex_\mathcal{P}(n, S_{3,5})$, and give the following theorem.

\textbf{Theorem 1.1.} Let $G$ be an $S_{3,5}$-free planar graph on $n$ vertices. Then $ex_\mathcal{P}(n, S_{3,5})$ $\leq$ $\tfrac{23n}{8}-\tfrac{9}{2}$ for all $n \geq 2$. The upper bound is tight when $n=12$.

The remainder of this paper is organized as follows. Essential preliminaries are provided in Section 2. Section 3 establishes the proof of Theorem 1.1, with detailed analysis of its implications.

\begin{figure}[h]
    \centering
    \includegraphics[width=0.45\textwidth]{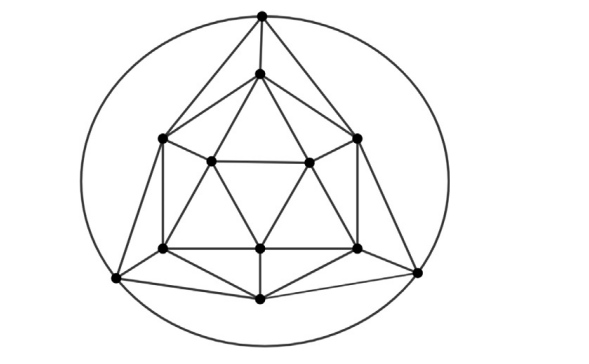}
    \caption{5-regular maximal planar graph on 12 vertices. }
    \label{Fig:11111} 
\end{figure}

\section{Preliminaries}

The upper bound in Theorem 1.1 is tight for 
$n=12$, as shown by the 5-regular maximal planar graph with 12 vertices and 30 edges in Figure 1. This graph contains no $S_{3,5}$, since every vertex has degree 5—whereas $S_{3,5}$ includes a vertex of degree 6.

\vspace{3mm}
\noindent\textbf{Claim 2.1. }If $G$ is a graph on $n$ vertices for $2\leq n\leq 12$, then the number of edges in $G$ is at most $\tfrac{23n}{8}-\tfrac{9}{2}$.

\noindent\textbf{Proof.} It is easy to see that a maximum planar graph with $n$ vertices contains $3n-6$ edges. Note that $3n-6$ $\leq$ $\tfrac{23n}{8}-\tfrac{9}{2}$ for $2\leq n\leq 12$. Thus, there is $e(G)$  $\leq$ $\tfrac{23n}{8}-\tfrac{9}{2}$ for $2\leq n\leq 12$. $\hfill\square$

\vspace{3mm}
For the convenience of discussion, according to Claim 2.1, we may further assume $n\geq 12$.

\vspace{3mm}
\noindent\textbf{Claim 2.2. } Let $G$ be an $S_{3,5}$-free graph. If there exist an edge cut set $E[S_1,S_2]$ in the graph $G$ such that $V(G)=S_1 \cup S_2 $, $e[S_1,S_2]$ $\leq$ $4$ and $|S_1|$, $|S_2| \geq 1$, then $e(G)$  $\leq$ $\tfrac{23n}{8}-\tfrac{9}{2}$.

\noindent\textbf{Proof.} Suppose that there exists an edge cut set $E[S_1,S_2]$ in the graph $G$, where $V(G)=S_1 \cup S_2 $, $e[S_1,S_2]$ $\leq$ $4$. Since $G$ is $S_{3,5}$-free, then $G[S_1]$ and $G[S_2]$ are $S_{3,5}$-free. 

By the induction hypothesis, $e(G[S_1])$ $\leq$ 
   $\tfrac{23|S_1|}{8}-\tfrac{9}{2}$, $e(G[S_2])$ $\leq$ $\tfrac{23|S_2|}{8}-\tfrac{9}{2}$. Then 	
\begin{align*}
e(G)&=e(G[S_1])+e(G[S_2])+e[S_1,S_2]\\
 &\leq \tfrac{23|S_1|}{8}-\tfrac{9}{2}+\tfrac{23|S_2|}{8}-\tfrac{9}{2}+4\\     
&\leq \tfrac{23n}{8}-5.
\end{align*}$\hfill\square$

\vspace{3mm}      
 In the following, assume that all edge cut sets $E[S_1,S_2]$ in the graph $G$ satisfy $e[S_1,S_2]$ $\geq$ $5$.

If $G$ is disconnected, there is an edge cut set $E[S_1,S_2]$ such that $e[S_1,S_2]$ = $0$ and $|S_1|$, $|S_2| \geq 1$. Then $e(G)$ $\leq$ $\tfrac{23n}{8}-\tfrac{9}{2}$ by the Claim 2.2. So assume that $G$ is connected. 

In all the following cases, let $G$ be an $S_{3,5}$-free planar graph with $n$ vertices for $n\geq 12$.

\vspace{3mm}
\noindent\textbf{Claim 2.3. }In the graph $G$, a vertex with degree at least 9 cannot be adjacent to a vertex with degree at least 4.

\noindent\textbf{Proof.} Suppose not. Let $xy$ be an edge in $G$ such that $d(x) \geq 9$ and $d(y) \geq 4$. Obviously, there are at least three vertices in $V (G)\backslash \{x\}$, say $y_1$, $y_2$ and $y_3$, which are adjacent to $y$. Because of $|N(x)\backslash y|$ $\geq$ 8, there
are at least five vertices $x_1$, $x_2$, $x_3$, $x_4$ and $x_5$, respectively, not in \{$y$, $y_1$, $y_2$, $y_3$\} which are adjacent to $x$. This implies we
got an $S_{3,5}$ in $G$ with backbone $xy$ and leaf-sets \{$x_1$, $x_2$, $x_3$, $x_4$, $x_5$\} and \{$y_1$, $y_2$, $y_3$\}, respectively, which
is a contradiction. This completes the proof.$\hfill\square$

\vspace{3mm}
\noindent\textbf{Claim 2.4. }If $\Delta(G)$$\geq$9 in the graph $G$, then  $e(G)$  $\leq$ $\tfrac{23n}{8}-\tfrac{9}{2}$.

\noindent\textbf{Proof.} Suppose that $d(x)$=$\Delta(G)$ with $d(x) \geq 9$. By Claim 2.3, we know that for any $y\in N(x)$ there is $d(y) \leq 3$. Without loss of generality, taking $x$ and nine neighbors of $x$ in $N(x)$ as a subset $H$. Deleting the vertices in $H$ from the graph $G$, the number of edges deleted is at most $9+9\times2=27$.
	
By the induction hypothesis, $e(G\backslash H)$ $\leq$ $\tfrac{23(n-10)}{8}-\tfrac{9}{2}$. Then
\begin{align*}
 e(G)&\leq e(G\backslash H)+27\\
 &\leq \tfrac{23(n-10)}{8}-\tfrac{9}{2}+27\\
&\leq \tfrac{23n}{8}-\tfrac{9}{2}-\tfrac{7}{4}.\end{align*}$\hfill\square$

\vspace{3mm}
\noindent\textbf{Claim 2.5. }If the graph $G$ contains a vertex $v$ with $d(v) = 8$, then $e(G)$  $\leq$ $\tfrac{23n}{8}-\tfrac{9}{2}$.

\noindent\textbf{Proof.}  Let $d(v) = 8$, $S_1$=\{$x_1$, $x_2$, $x_3$, $x_4$, $x_5$, $x_6$, $x_7$, $x_8$\} as the set of its eight neighbors, and define $H=S_1 \cup \{v\} $. If a vertex in $S_1$ is adjacent to at least two vertices within $S_1$, then it has no neighbors in $V(G)\backslash H$. If a vertex is adjacent to one vertex within $S_1$, then it has at most one neighbor in $V(G)\backslash H$. If a vertex is not adjacent to any vertices within $S_1$, then it has at most two neighbors in $V(G)\backslash H$. Next, we conduct a classification discussion on the total number of vertices in $S_1$ that are adjacent to at least two vertices within $S_1$.

\textbf{Case 1.} If each vertex in $S_1$ has at most one neighbor within $S_1$, then removing $H$ from the graph $G$ deletes at most $8+8\times2$=$24$ edges.

\textbf{Case 2.} If there is one vertex in $S_1$ that has at least two neighbors within $S_1$, and the remaining 7 vertices have at most one neighbor within $S_1$, then removing $H$ from the graph $G$ deletes at most $8+2+2+5\times2$=$22$ edges.

\textbf{Case 3.} If there is more than one vertex in $S_1$ that has at least two neighbors within $S_1$, and the remaining vertices have at most one neighbor within $S_1$, then removing $H$ from the graph $G$ deletes at most $3(m+1)-6+(8-m)+(8-m)\times2$=21 edges. 

At this case, let $V_m$ be the set containing these vertices. Therefore, $V_m\cup\{v\}$ contains at least three vertices. In this case, we consider the maximum number of edges in the subgraph induced by $V_m\cup\{v\}$. A planar graph with $n \geq 3$ vertices has at most $3n-6$ edges. Let $m$ be the number of vertices in $V_m$. We have $e(V_m\cup\{v\})$$\leq$ $3(m+1)-6$. Additionally, beyond the edges mentioned earlier, there are at most $(8-m)+(8-m)\times2$ extra edges.

To conclude, deleting $H$ without considering the condition that the size of the edge cut set is at least 5 will delete at most 24 edges. Thus, under the condition that the size of the edge cut set is at least 5,  the number of edges deleted is at most 24.  

By the induction hypothesis, $e(G\backslash H)$ $\leq$  $\tfrac{23(n-9)}{8}-3$,
\begin{align*}
e(G)&\leq e(G\backslash H)+24\\ 
	&\leq \tfrac{23(n-9)}{8}-\tfrac{9}{2}+24\\ 
	&\leq \tfrac{23n}{8}-\tfrac{9}{2} -\tfrac{15}{8}.
\end{align*}
$\hfill\square$

\vspace{3mm}
\noindent\textbf{Claim 2.6. } Let $v$ be a vertex with $d(v) = 7$ in the graph $G$ and the set consisting of the seven neighbors of $v$ be $S_1$=\{$x_1$, $x_2$, $x_3$, $x_4$, $x_5$, $x_6$, $x_7$\}. If at least one vertex in $S_1$ has no neighbors within $S_1$, then $e(G)$  $\leq$ $\tfrac{23n}{8}-\tfrac{9}{2}$.

\noindent\textbf{Proof:} Let $d(v) = 7$ and define $H=S_1 \cup \{v\} $. If a vertex in $S_1$ is adjacent to at least one vertex within $S_1$, then it has at most one neighbor in $V(G)\backslash H$. If a vertex is not adjacent to any vertices within $S_1$, then it has at most two neighbors in $V(G)\backslash H$. Next, we conduct a classification discussion on the total number of vertices in $S_1$ that are adjacent to at least one vertex within $S_1$.

\textbf{Case 1.} If each vertex in $S_1$ has no neighbors within $S_1$, then removing $H$ from the graph $G$ deletes at most $7+7\times2$=$21$ edges.

\textbf{Case 2.} If there is one vertex in $S_1$ that have at least one neighbor within $S_1$, and the remaining six vertices have no neighbors within $S_1$, then it is impossible. 

\textbf{Case 3.} If there are two vertices in $S_1$ that have at least one neighbor within $S_1$, and the remaining five vertices have no neighbors within $S_1$, then removing $H$ from the graph $G$ deletes at most $7+1+5\times2+2\times1$=$20$ edges.

\textbf{Case 4.} If there are three vertices in $S_1$ that have at least one neighbor within $S_1$, and the remaining four vertices have no neighbors within $S_1$, then removing $H$ from the graph $G$ deletes at most $7+3+4\times2+3\times1$=$21$ edges. 

\textbf{Case 5.} If there are four vertices in $S_1$ that have at least one neighbor within $S_1$, and the remaining three vertices have no neighbors within $S_1$, then removing $H$ from the graph $G$ deletes at most $9+3+3\times2+4\times1$=$22$ edges.

In this case, these four vertices and vertex $v$ form a planar graph with 5 vertices. 
For a planar graph with 5 vertices, the maximum number of edges is $3\times5-6$=$9$. To preserve planarity during the removal of a subgraph 
$H$, at most $3+3\times2+4\times1$=$13$  additional edges beyond these 9 edges must be deleted.

\textbf{Case 6.} If there are five vertices in $S_1$ that have at least one neighbor within $S_1$, and the remaining two vertices have no neighbors within $S_1$, then removing $H$ from the graph $G$ deletes at most $12+2+2\times2+5\times1$=$23$ edges.

In this case, these five vertices and vertex $v$ form a planar graph with 6 vertices.
For a planar graph with 6 vertices, the maximum number of edges is $3\times6-6$=$12$. To preserve planarity during the removal of a subgraph 
$H$, at most $2+2\times2+5\times1$=$11$   additional edges beyond these 12 edges must be deleted.

\textbf{Case 7.} If there are six vertices in $S_1$ that have at least one neighbor within $S_1$, and the remaining one vertex have no neighbors within $S_1$, then removing $H$ from the graph $G$ deletes at most $15+1+1\times2+6\times1$=$24$ edges. 

Without loss of generality, let these six vertices be \{$x_1$, $x_2$, $x_3$, $x_4$, $x_5$, $x_6$\}.
In this case, \{$x_1$, $x_2$, $x_3$, $x_4$, $x_5$, $x_6$\}$\cup \{v\}$ form a planar graph with 7 vertices. For a planar graph on 7 vertices achieving the maximum of 15 edges ($3\times7-6$=$15$) with a vertex $v$ adjacent to all six other vertices, the graph is necessarily isomorphic to that shown in Figure 2.
\begin{figure}[h]
    \centering
    \includegraphics[width=0.45\textwidth]{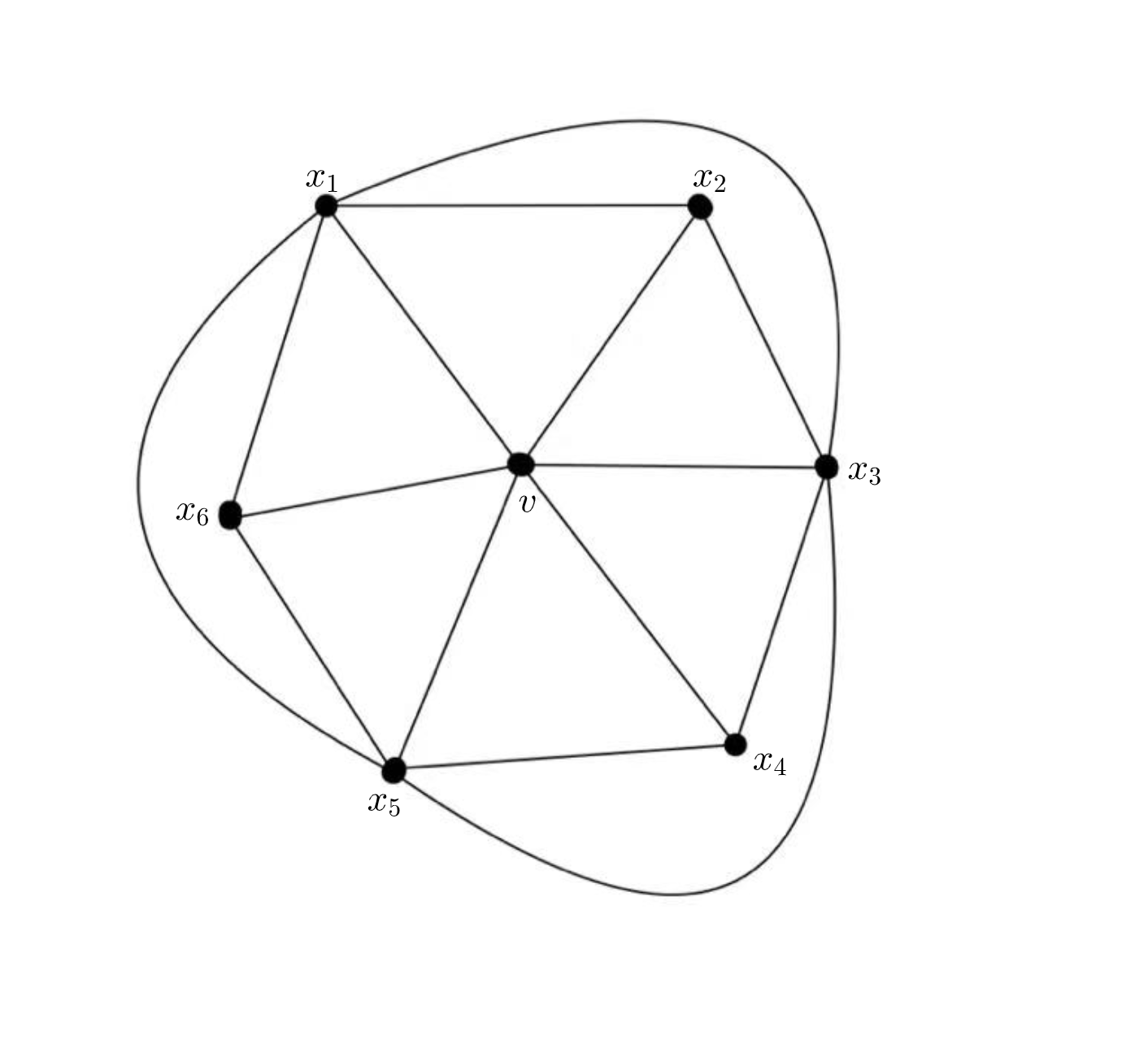}
    \caption{A planar graph with 7 vertices, including a central vertex $v$, has at most 15 edges. }
    \label{Fig:11111} 
\end{figure}

In this case, regardless of which face the vertex $x_7$ lies in, there exists an edge cut set with size at most 2. Therefore, the number of edges deleted is at most 23.

To conclude, if at least one vertex in $S_1$ has no neighbors within $S_1$, then removing $H$ from the graph $G$ deletes at most 23 edges.

By the induction hypothesis, 
$e(G\backslash H)$ $\leq$  $\tfrac{23(n-8)}{8}-\tfrac{9}{2}$,
    
    $e(G)$$\leq$ $e(G\backslash H)$+$23$ 
    $\leq$  $\tfrac{23(n-8)}{8}-\tfrac{9}{2} +23$ 
    $\leq$  $\tfrac{23n}{8}-\tfrac{9}{2}$.$\hfill\square$

\vspace{3mm}
\noindent\textbf{Claim 2.7. } Let $v$ be a vertex with $d(v) = 6$ in the graph $G$, then $e(G)$  $\leq$ $\tfrac{23n}{8}-\tfrac{9}{2}$.

\noindent\textbf{Proof:} Let $d(v) = 6$ and define the set consisting of the six neighbors of $v$ be $S_1$=\{$x_1$, $x_2$, $x_3$, $x_4$, $x_5$, $x_6$\}. $H=S_1 \cup \{v\}$. Regardless of the adjacency within $S_1$, each vertex in $S_1$ has at most two neighbors in $V(G)\backslash H$. Then removing $H$ from the graph $G$ deletes at most $3\times 7-6$+$6\times 2$=27 edges.

In this case, $H$ form a planar graph with 7 vertices. For a planar graph with 7 vertices, the maximum number of edges is $3\times7-6$=$15$. To preserve planarity during the removal of a subgraph 
$H$, at most $6\times2$=$12$ additional edges beyond these 15 edges must be deleted. Similarly, for a planar graph on 7 vertices achieving the maximum of 15 edges ($3\times7-6$=$15$) with a vertex $v$ adjacent to all six other vertices, the graph is necessarily isomorphic to that shown in Figure 2.

At this time, the vertices $\{x_1, x_3, x_5\}$ form one equivalence class, while $\{x_2, x_4, x_6\}$ form another. Similarly, the edge sets $\{x_1x_2, x_2x_3, x_3x_4, x_4x_5, x_5x_6\}$ and $\{x_1x_3, x_3x_5, x_1x_5\}$ each constitute distinct equivalence classes.

Based on the properties of isomorphism, we only need to consider the situation of Figure 2. 

\noindent\textbf{Claim:} Note that if the vertices of $\{x_2, x_4, x_6\}$ lie in the distinct faces and each has at least one neighbor in $V(G)\backslash H$, then there exist a cut set of size at most 4, a contraction. 

\noindent\textbf{Case 1:} For Figure 2, if we delete 27 edges, each vertex of $H$ has two neighbors of $V(G)\backslash H$, there must exist a cut set of size at most 4,a contraction.

\noindent\textbf{Case 2:} If we delete 26 edges, we need only consider removing one additional edge in Case 1 since Case 1 contains exactly 27 edges satisfying $e(G[H]) + e(H, V(G)\setminus H) = 27$. This leads to three subcases:

\textbf{Case 2.1:} If we delete a neighbor of a vertex in
$S_1$ within $V(G)\backslash H$ – regardless of which neighbor in $V(G)\backslash H$ is removed - then the vertices $\{x_2, x_4, x_6\}$ lie in distinct faces and each has at least one neighbor in $V(G)\backslash H$. By claim, this implies the existence of an edge cut set of size at most 4, a contradiction.

\textbf{Case 2.2:} Similarly, if we remove any single edge from $\{x_1x_2, x_2x_3, x_3x_4, x_4x_5, x_5x_6\}$, then by Claim, this implies the existence of an edge cut set of size at most 4, yielding a contradiction.

\textbf{Case 2.3:} Also, if we remove any single edge from $\{x_1x_3, x_3x_5, x_1x_5\}$, then by Claim, this implies the existence of an edge cut set of size at most 4, yielding a contradiction.

\noindent\textbf{Case 3:} Similarly, if we remove 25 edges, we need only consider removing two additional edges in Case 1. Therefore, we have four subcases:

\textbf{Case 3.1:} If we delete a neighbor of a vertex in
$\{x_1, x_3, x_5\}$ within $V(G)\backslash H$ at first. Assume that we delete a neighbor of vertex $x_1$ within $V(G)\backslash H$.

\textbf{Case 3.1.1:} If we delete two neighbors of vertex $x_1$ within $V(G)\backslash H$. By claim, this implies the existence of an edge cut set of size at most 4, yielding a contradiction.

\textbf{Case 3.1.2:} Moreover, we delete one neighbor of a vertex in $S_1 \backslash$ $\{x_1\}$ within $V(G)\backslash H$. By claim, this implies the existence of an edge cut set of size at most 4, yielding a contradiction.

\textbf{Case 3.1.3:} Moreover, we delete one edge of $\{x_1x_2, x_2x_3, x_3x_4, x_4x_5, x_5x_6\}$ $\cup$ $\{x_1x_3, x_3x_5, x_1x_5\}$ . By Claim, this implies the existence of an edge cut set of size at most 4, producing a contradiction.

\textbf{Case 3.2:} If we delete a neighbor of a vertex in
$\{x_2, x_4, x_6\}$ within $V(G)\backslash H$ at first. Assume that we delete a neighbor of vertex $x_2$ within $V(G)\backslash H$.

\textbf{Case 3.2.1:} If we delete two neighbors of vertex $x_2$ within $V(G)\backslash H$. Then there must exist two distinct components $H_1$ in the face $x_1x_5x_6$ and $H_2$ in the face $x_3x_4x_5$.

Since $|V(G)|$$\geq$$12$, we just consider the following two situations:

(i) If $|V(H_1)|$=1 and $|V(H_2)|$$\geq$$2$ (or $|V(H_1)|$=1 and $|V(H_1)|$$\geq$$2$), then by the induction hypothesis:

$e(G\backslash(H \cup H_1\cup  H_2 ))$ $\leq$ $\tfrac{23(n-8-|V(H_1)|)}{8}-\tfrac{9}{2}$, $e(H_2)$$\leq$$\tfrac{23|V(H_2)|}{8}-\tfrac{9}{2}$.
\vspace{1mm}

Then $e(G)$$\leq$$\tfrac{23(n-8-|V(H_1)|)}{8}-\tfrac{9}{2}$+$\tfrac{23|V(H_1)|}{8}-\tfrac{9}{2}$+25
	$\leq$$\tfrac{23n}{8}-\tfrac{9}{2}$.
\vspace{1mm}

(ii) If $|V(H_1)|$$\geq$$2$ and $|V(H_2)|$$\geq$$2$, then by the induction hypothesis:
\vspace{1mm}

$e(G\backslash(H \cup H_1\cup  H_2 ))$ $\leq$ $\tfrac{23(n-7-|V(H_1)|-|V(H_2)|)}{8}-\tfrac{9}{2}$, $e(H_1)$$\leq$$\tfrac{23|V(H_1)|}{8}-\tfrac{9}{2}$, $e(H_2)$$\leq$$\tfrac{23|V(H_2)|}{8}-\tfrac{9}{2}$.
\vspace{1mm}

Then $e(G)$$\leq$ $\tfrac{23(n-7-|V(H_1)|-|V(H_2)|)}{8}-\tfrac{9}{2}$+$\tfrac{23|V(H_1)|}{8}-\tfrac{9}{2}$+$\tfrac{23|V(H_2)|}{8}-\tfrac{9}{2}$+25
	$\leq$$\tfrac{23n}{8}-\tfrac{9}{2}$.
\vspace{1mm}

\textbf{Case 3.2.2:} Moreover, we delete one neighbor of a vertex in $S_1 \backslash$ $\{x_2\}$ within $V(G)\backslash H$. By claim, this implies the existence of an edge cut set of size at most 4, yielding a contradiction.

\textbf{Case 3.2.3:} Moreover, we delete one edge of $\{x_1x_2, x_2x_3, x_3x_4, x_4x_5, x_5x_6\}$ $\cup$ $\{x_1x_3, x_3x_5, x_1x_5\}$. By Claim, this implies the existence of an edge cut set of size at most 4, producing a contradiction.

\textbf{Case 3.3:} If we delete one edge of 
$\{x_1x_2, x_2x_3, x_3x_4, x_4x_5, x_5x_6\}$ at first. Assume that we delete the edge  $x_1x_2$.

\textbf{Case 3.3.1:} Moreover, we delete one neighbor of a vertex in $S_1$ within $V(G)\backslash H$. By claim, this implies the existence of an edge cut set of size at most 4, yielding a contradiction.

\textbf{Case 3.3.2:} Moreover, we delete one edge of $\{x_2x_3, x_3x_4, x_4x_5, x_5x_6\}$ $\cup$ $\{x_1x_3, x_3x_5, x_1x_5\}$. By Claim, this implies the existence of an edge cut set of size at most 4, producing a contradiction.

\textbf{Case 3.4:} If we delete one edge of 
$\{x_1x_3, x_3x_5, x_1x_5\}$ at first. Assume that we delete the edge $x_1x_3$.

\textbf{Case 3.4.1:} Moreover, we delete one neighbor of a vertex in $S_1$ within $V(G)\backslash H$. By claim, this implies the existence of an edge cut set of size at most 4, yielding a contradiction.

\textbf{Case 3.4.2:} Moreover, we delete one edge of $\{x_1x_2,x_2x_3, x_3x_4, x_4x_5, x_5x_6\}$. By Claim, this implies the existence of an edge cut set of size at most 4, producing a contradiction.

\textbf{Case 3.4.3:} Moreover, we delete one edge of $\{x_3x_5, x_1x_5\}$. Assume that we delete the edge $x_1x_5$. Then there must exist two distinct components $H_1$ in the face $x_1x_2x_3x_5x_6$ and $H_2$ in the face $x_3x_4x_5$.

Since $|V(G)|$$\geq$$12$, we just consider the following two situations:
\vspace{1mm}

(i) If $|V(H_1)|$=1 and $|V(H_2)|$$\geq$$2$ (or $|V(H_1)|$=1 and $|V(H_1)|$$\geq$$2$), then by the induction hypothesis:
\vspace{1mm}

$e(G\backslash(H \cup H_1\cup  H_2 ))$ $\leq$ $\tfrac{23(n-8-|V(H_1)|)}{8}-\tfrac{9}{2}$, $e(H_2)$$\leq$$\tfrac{23|V(H_2)|}{8}-\tfrac{9}{2}$.
\vspace{1mm}

Then $e(G)$$\leq$$\tfrac{23(n-8-|V(H_1)|)}{8}-\tfrac{9}{2}$+$\tfrac{23|V(H_1)|}{8}-\tfrac{9}{2}$+25
	$\leq$$\tfrac{23n}{8}-\tfrac{9}{2}$.
\vspace{1mm}

(ii) If $|V(H_1)|$$\geq$$2$ and $|V(H_2)|$$\geq$$2$, then by the induction hypothesis:
\vspace{1mm}

$e(G\backslash(H \cup H_1\cup  H_2 ))$ $\leq$ $\tfrac{23(n-7-|V(H_1)|-|V(H_2)|)}{8}-\tfrac{9}{2}$, $e(H_1)$$\leq$$\tfrac{23|V(H_1)|}{8}-\tfrac{9}{2}$, $e(H_2)$$\leq$$\tfrac{23|V(H_2)|}{8}-\tfrac{9}{2}$.
\vspace{1mm}

Then $e(G)$$\leq$ $\tfrac{23(n-7-|V(H_1)|-|V(H_2)|)}{8}-\tfrac{9}{2}$+$\tfrac{23|V(H_1)|}{8}-\tfrac{9}{2}$+$\tfrac{23|V(H_2)|}{8}-\tfrac{9}{2}$+25
	$\leq$$\tfrac{23n}{8}-\tfrac{9}{2}$.
\vspace{1mm}

\noindent\textbf{Case 4:} Similarly, when removing between 16 and 24 edges, at least one but at most two components become adjacent to the subgraph $G[H]$. This follows because at least one vertex in $S_1$ has neighbors in $V(G)\backslash H$. If there are more than two distinct components, then $e[H,V(G)\backslash H]$ $\geq 15$ must hold to avoid an edge cut set of size at most 4, yielding a contradiction. We restrict our analysis to the following two cases:

\textbf{Case 4.1:} If just exists one component $H_1$ adjacents to the subgraph $G[H]$. Then $|V(H_1)|$$\geq$$2$, by the induction hypothesis:
\vspace{1mm}

$e(G\backslash(H \cup H_1))$ $\leq$ $\tfrac{23(n-7-|V(H_1)|)}{8}-\tfrac{9}{2}$, $e(H_1)$$\leq$$\tfrac{23|V(H_1)|}{8}-\tfrac{9}{2}$.
\vspace{1mm}

Then $e(G)$$\leq$$\tfrac{23(n-7-|V(H_1)|)}{8}-\tfrac{9}{2}$+$\tfrac{23|V(H_1)|}{8}-\tfrac{9}{2}$+24
	$\leq$$\tfrac{23n}{8}-\tfrac{9}{2}$.
\vspace{1mm}

\textbf{Case 4.2:} If just exist two components $H_1$ and $H_2$ adjacent to the subgraph $G[H]$.

Since $|V(G)|$$\geq$$12$, we just consider the following two situations:
\vspace{1mm}

(i) If $|V(H_1)|$=1 and $|V(H_2)|$$\geq$$2$ (or $|V(H_1)|$=1 and $|V(H_1)|$$\geq$$2$), then by the induction hypothesis:
\vspace{1mm}

$e(G\backslash(H \cup H_1\cup  H_2 ))$ $\leq$ $\tfrac{23(n-8-|V(H_1)|)}{8}-\tfrac{9}{2}$, $e(H_2)$$\leq$$\tfrac{23|V(H_2)|}{8}-\tfrac{9}{2}$.
\vspace{1mm}

Then $e(G)$$\leq$$\tfrac{23(n-8-|V(H_1)|)}{8}-\tfrac{9}{2}$+$\tfrac{23|V(H_1)|}{8}-\tfrac{9}{2}$+24
	$\leq$$\tfrac{23n}{8}-\tfrac{9}{2}$.
\vspace{1mm}

(ii) If $|V(H_1)|$$\geq$$2$ and $|V(H_2)|$$\geq$$2$, then by the induction hypothesis:
\vspace{1mm}

$e(G\backslash(H \cup H_1\cup  H_2 ))$ $\leq$ $\tfrac{23(n-7-|V(H_1)|-|V(H_2)|)}{8}-\tfrac{9}{2}$, $e(H_1)$$\leq$$\tfrac{23|V(H_1)|}{8}-\tfrac{9}{2}$, $e(H_2)$$\leq$$\tfrac{23|V(H_2)|}{8}-\tfrac{9}{2}$.
\vspace{1mm}

Then $e(G)$$\leq$ $\tfrac{23(n-7-|V(H_1)|-|V(H_2)|)}{8}-\tfrac{9}{2}$+$\tfrac{23|V(H_1)|}{8}-\tfrac{9}{2}$+$\tfrac{23|V(H_2)|}{8}-\tfrac{9}{2}$+25
	$\leq$$\tfrac{23n}{8}-\tfrac{9}{2}$.
\vspace{1mm}

\noindent\textbf{Case 5:} Obviously, when removing at most 15 edges, by the induction hypothesis, Claim 2.7 is ture.

\vspace{3mm}
\noindent\textbf{Claim 2.8. } If the graph $G$ contains a 4-7 edge, then $e(G)$  $\leq$ $\tfrac{23n}{8}-\tfrac{9}{2}$.

\noindent\textbf{Proof.} Let $xy$$\in$$E(G)$ be a 7-4 edge. There are at least 2 triangles containing the edge $xy$. Otherwise, $G$ contains an $S_{3,5}$. We subdivide the cases based on the number of triangles containing the edge $xy$.

\textbf{Case 1. There are 3 triangles containing the edge $xy$.} 

Let $a$, $b$ and $c$ be the vertices in $G$ that are adjacent to both $x$ and $y$. Let $d$, $e$ and $f$ be the
vertices adjacent to $x$ but not to $y$, as shown in the Figure 3. Define $S_1$=\{$a$, $b$, $c$\}, $S_2$=\{$d$, $e$, $f$\}, $S_3$=$S_1\cup S_2\cup \{y\}$ and let $H$=$S_3\cup \{x\}$. 

\begin{figure}[h]
    \centering
    \includegraphics[width=0.5\textwidth]{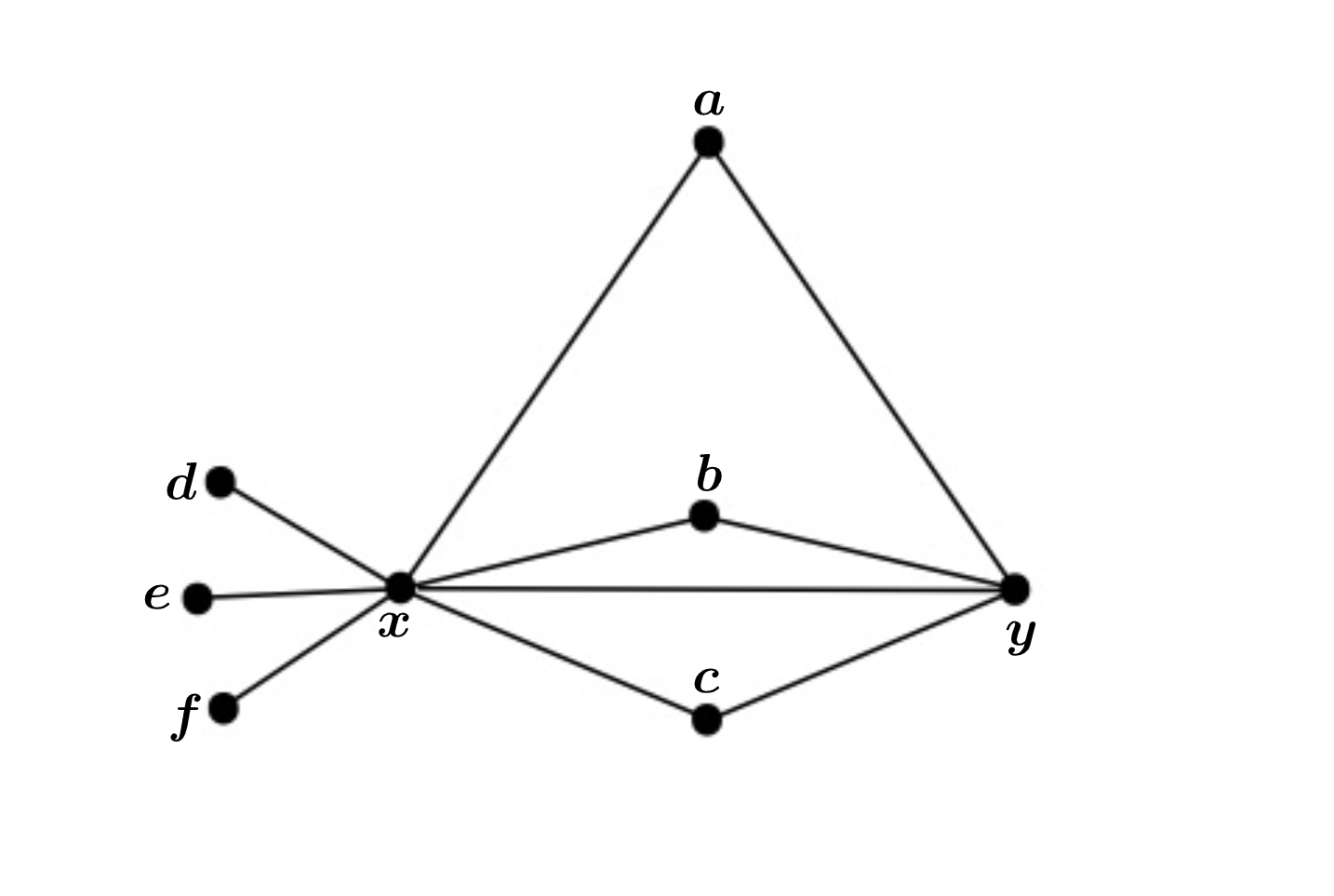}
    \caption{The 7-4 edge $xy$ is contained in 3 triangles.}
    \label{Fig:57354} 
\end{figure}

Removing the vertices in $H$ from the graph $G$. Assume that $a$ has two neighbors in $V(G)\backslash H$, we immediately get an $S_{3,5}$. Thus, any vertices in $S_1$ can have at most one neighbor in $V(G)\backslash H$.  
It can form a path of length 2 within $S_1$ and the edges of $S_1$ do not affect the neighbors of the vertices of $S_1$ in $V(G)\backslash H$. Each of the vertices $d$, $e$ and $f$ can have at most two neighbors in $V(G)\backslash H$. If a vertex $s$ in $S_2$ is adjacent to a vertex $t$ in $S_1$, this adjacent situation does not alter the set of vertices adjacent to $t$ within $S_1$ in the subgraph $V(G)\backslash H$. If two of the vertices of $\{d, e, f\}$ are adjacent, this implies that each of them has at most one neighbor in $V(G)\backslash H$. If $d$ is adjacent to at least one vertex in $S_1 \cup \{e, f\}$, then it has at most one neighbor in $V(G)\backslash H$. Similarly, the same situation hold for $e$ and $f$.

Here, $x$ is a vertex of degree 7 in $G$. According to Claim 2.6, it is only necessary to consider the case where all the vertices in $S_3$ have at least one other neighbor within $S_3$. Specifically, consider vertices $d$, $e$ and $f$ in $S_3$.

When removing the vertices of $H$ from $G$, at most $10+2+3\times1+(1+2)\times3$=24 edges are deleted. Since a planar graph on eight vertices has at most $3\times8-6$=18 edges and the set $S_1\cup S_2$ contains vertices with at most one neighbor in $V(G)\backslash H$, there exist at most six additional edges beyond these 18 edges.
At this time, $d$, $e$ and $f$ are adjacent to at least two other vertices in $S_1\cup S_2$, respectively, the three vertices $d$, $e$ and $f$ are connected to at most 9 edges in total. If we delete 24 edges, then all the vertices in $S_1\cup S_2$ have one neighbor in $V(G)\backslash H$. The possible cases as shown in the Figure 4.

\begin{figure}[h]
    \centering
    \includegraphics[width=0.9\textwidth]{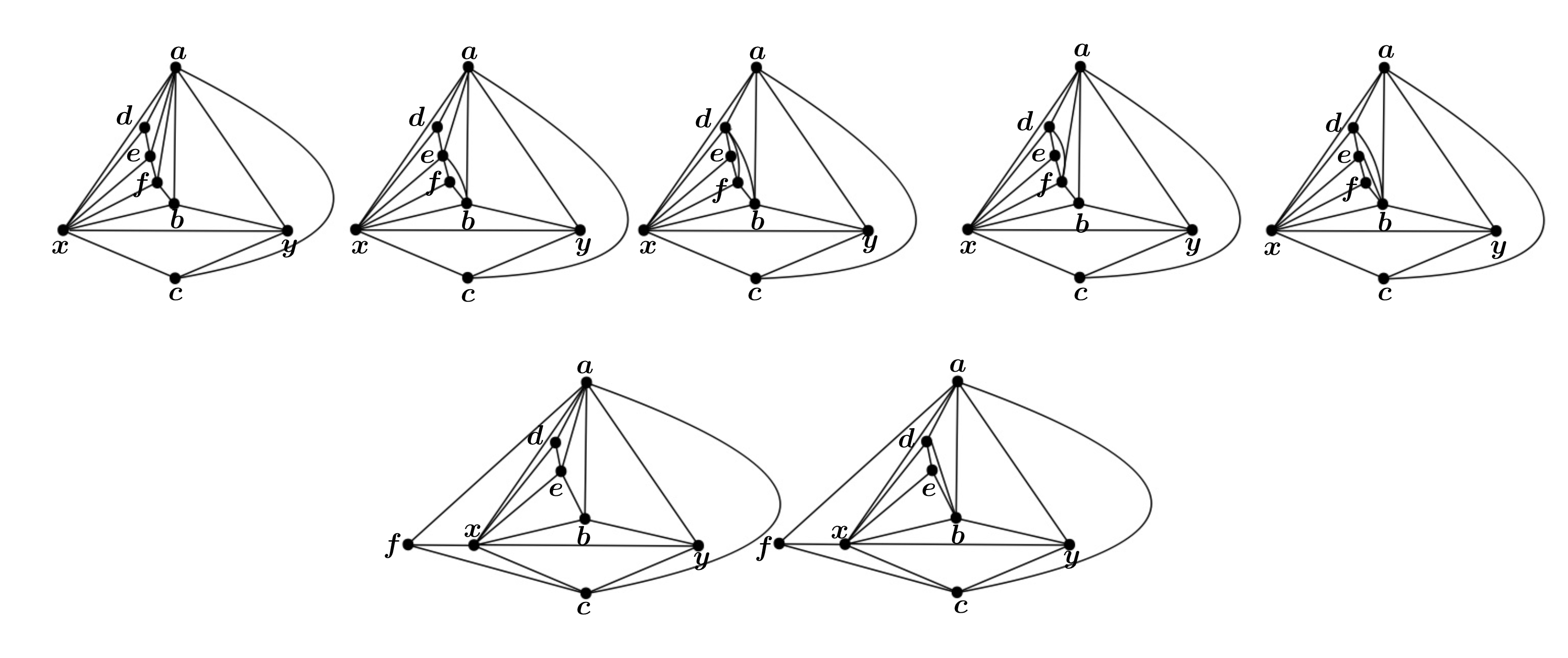}
    \caption{If we delete 24 edges, here are the adjacent situations of the vertices of $H$.}
    \label{Fig:57444} 
\end{figure}

For Figure 4, if we remove 24 edges, there must exist an edge cut set with size at most 3. Therefore, we must delete at most 23 edges to avoid finding an edge cut set with size at most 4.

By the induction hypothesis, Claim 2.8 holds. 

\textbf{Case 2. There are 2 triangles containing the edge $xy$.} 

Let $a$ and $b$ be the vertices in $G$ that are adjacent to both $x$ and $y$. Let $c$, $d$, $e$ and $f$ be the
vertices adjacent to $x$ but not to $y$, and $m$ be the
vertex adjacent to $y$ but not to $x$, as shown in the Figure 5. Define $S_1$=\{$a$, $b$\}, $S_2$=\{$c$, $d$, $e$, $f$\}, $S_3$=$S_1\cup S_2\cup\{y\}$ and let $H$=$S_3\cup\{x,m\}$. 

\begin{figure}[h]
    \centering
    \includegraphics[width=0.5\textwidth]{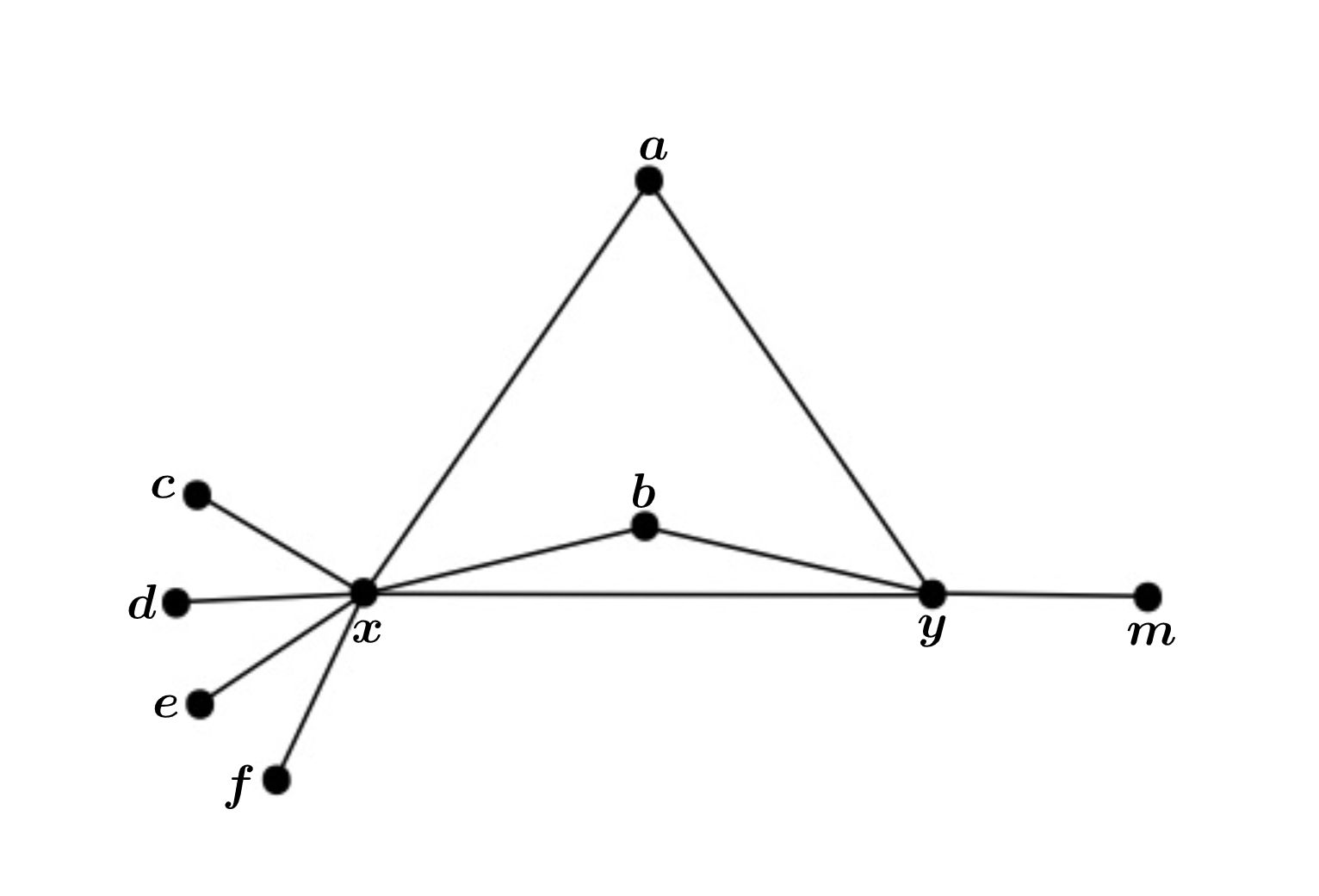}
    \caption{The 7-4 edge $xy$ is contained in 2 triangles.}
    \label{Fig:57353} 
\end{figure}

Removing the vertices in $H$ from the graph $G$. Assume that $a$ has two neighbors in $V(G)\backslash H$, we immediately get an $S_{3,5}$. Thus, any vertices in $S_1$ can have at most one neighbor in $V(G)\backslash H$.  
It can form a path of length 1 within $S_1$ and the edge of $S_1$ do not affect the neighbors of the vertices of $S_1$ in $V(G)\backslash H$. Each vertex in $S_2$ can have at most two neighbors in $V(G)\backslash H$.
If a vertex $r$ in $S_2$ is adjacent to a vertex $p$ in $S_1$, this adjacent situation does not alter the set of vertices adjacent to $p$ within $S_1$ in the subgraph $V(G)\backslash H$. If a vertex $r$ in $S_2$ is adjacent to at least one vertex in $S_1 \cup S_2 $, it makes these vertices in $S_2$ have at most one neighbor in $V(G)\backslash H$. The vertex $m$ can have at most four neighbors in $V(G)\backslash H$. If $m$ is adjacent to a vertex in $S_1$, then this vertex in $S_1$ have no neighbors in $V(G)\backslash H$, and $m$ has at most four neighbors in $V(G)\backslash H$. If $m$ is adjacent to a vertex in $S_2$, then this vertex in $S_2$ has at most one neighbor in $V(G)\backslash H$, and $m$ has at most three neighbors in $V(G)\backslash H$. Therefore, when we delete as many edges as possible, assume that $m$ is not adjacent to any vertices of $S_2$.

Here, $x$ is a vertex of degree 7. According to Claim 2.6, it is only necessary to consider the case where all the vertices in $S_3$ have at least one neighbor within $S_3$. Specifically, consider vertices $c$, $d$, $e$ and $f$ in $S_3$.

Deleting the vertices in $H$ from the graph $G$, the number of edges deleted is at most $10+1+2\times1+(1+2)\times4+4$=29.

For any planar embedding of the vertex set $S_3\cup x$ (with $|S_3\cup x|=8$), the edge count is bounded by $3\times8-6$=18 due to Euler's formula. Given that each vertex in $S_1\cup S_2$ has at most one neighbor in $V(G)\backslash H$, there are at most six additional edges beyond these 18 edges. Furthermore, vertex $m$ has at most five incident edges. Therefore, the maximum number of edges removed is bounded by 29. At this case, $c$, $d$, $e$, and $f$ are adjacent to at least two other vertices in $S_1\cup S_2$, respectively. The vertices $c$, $d$, $e$, and $f$ are in the same face and these vertices connected to at most 12 edges in total. The possible situations as shown in the Figure 6:

\begin{figure}[h]
    \centering
   \includegraphics[width=0.9\textwidth]{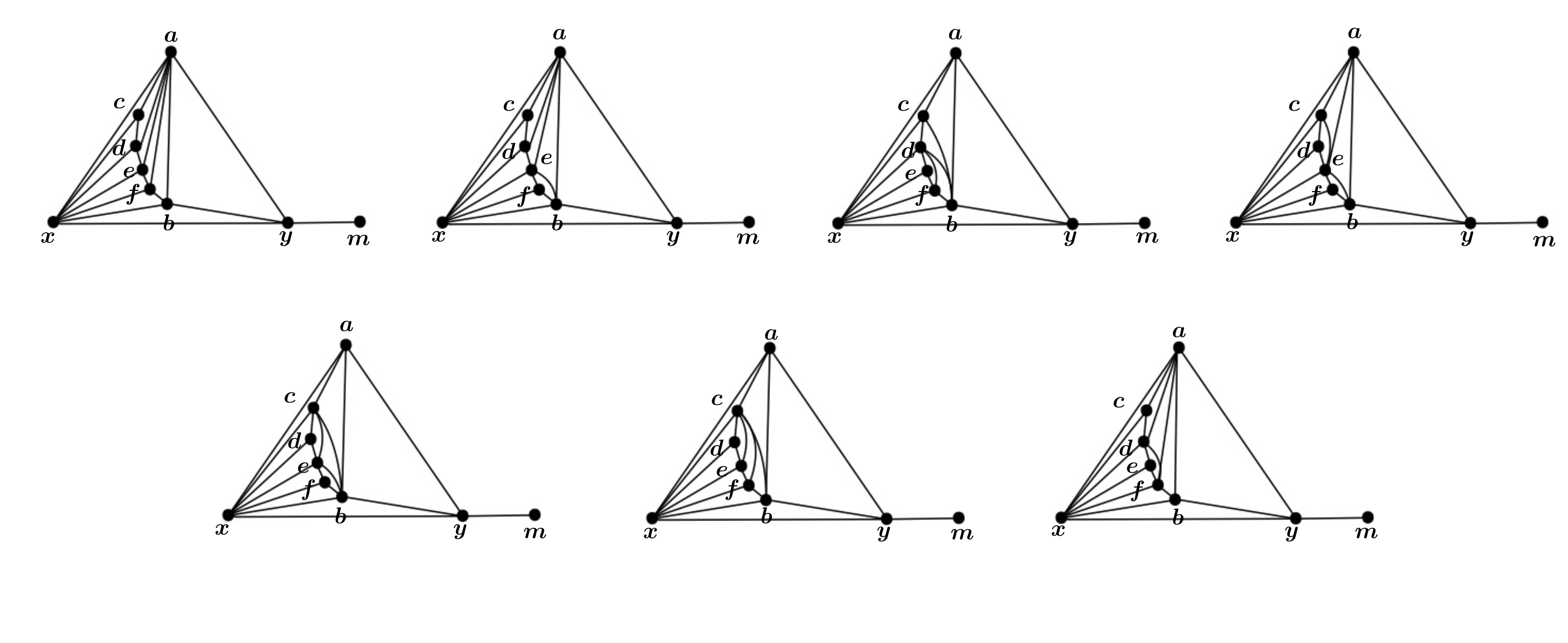}
   \caption{If we delete 29 edges, here are the adjacent situations of the vertices of H.
}
    \label{Fig:5735} 
\end{figure}

When deleting 29 edges, the subgraph $G[S_3\cup {x}]$ is a triangulation which all faces have degree 3, so there must exist the edge cut set of size at most 3, a contraction. Therefore, the number of edges deletes is at most 28. When removing between 11 edges and 28 edges, there must exist a connected component $H_l$ in the resulting graph $G\backslash G[H] $, because there exists at least one vertex in $S_1\cup S_2\cup {m}$ that has at least one neighbor in $V(G)\backslash H$.

Suppose that there are $l$ vertices in $H_l$ for $l\geq1$. If $l=1$, then we can  removing the vertices of $H\cup H_l$, there are 10 vertices and 28 edges. By the induction hypothesis, Claim 2.8 holds.    

Otherwise $l\geq2$, $H_l$ is a subgraph of the graph $G$. Therefore, $H_l$ does not contain any $S_{3,5}$ as a subgraph. Because the number of vertices in $H_l$ is less than $n$. By the induction hypothesis, $H_l$ has at most $\tfrac{23l}{8}-\tfrac{9}{2}$ edges. Let $V_l$ be the vertex set of $H_l$. Then, if we delete $H \cup V_l $ from the graph $G$, we will delete at most $\tfrac{23l}{8}-\tfrac{9}{2}$+28 edges.

By the induction hypothesis,

$e(G\backslash(H \cup V_l))$ $\leq$ $\tfrac{23(n-(l+9))}{8}-\tfrac{9}{2}$+$\tfrac{23l}{8}-\tfrac{9}{2}$+28
	$\leq$$\tfrac{23n}{8}-\tfrac{9}{2}$.

Therefore, deleting at most 28 edges, Claim 2.8 holds. $\hfill\square$

\vspace{3mm}

\vspace{3mm}
\noindent\textbf{Claim 2.9. } If the graph $G$ contains a 5-7 edge, then $e(G)$  $\leq$ $\tfrac{23n}{8}-\tfrac{9}{2}$.

\noindent\textbf{Proof.} Let $xy$$\in$$E(G)$ be a 5-7 edge. There are at least 3 triangles containing the edge $xy$. Otherwise, $G$ contains an $S_{3,5}$. We subdivide the cases based on the number of triangles containing the edge $xy$.

\textbf{Case 1. There are 4 triangles containing the edge $xy$.} 

Let $a$, $b$, $c$ and $d$ be the vertices in $G$ that are adjacent to both $x$ and $y$. Let $e$ and $f$ be the
vertices adjacent to $y$ but not to $x$, as shown in the Figure 7. Define $S_1$=\{$a$, $b$, $c$, $d$\}, $S_2$=\{$e$, $f$\}, $S_3$=$S_1\cup S_2\cup \{x\}$ and let $H$=$S_1\cup \{y\}$. 

\begin{figure}[h]
    \centering
    \includegraphics[width=0.3\textwidth]{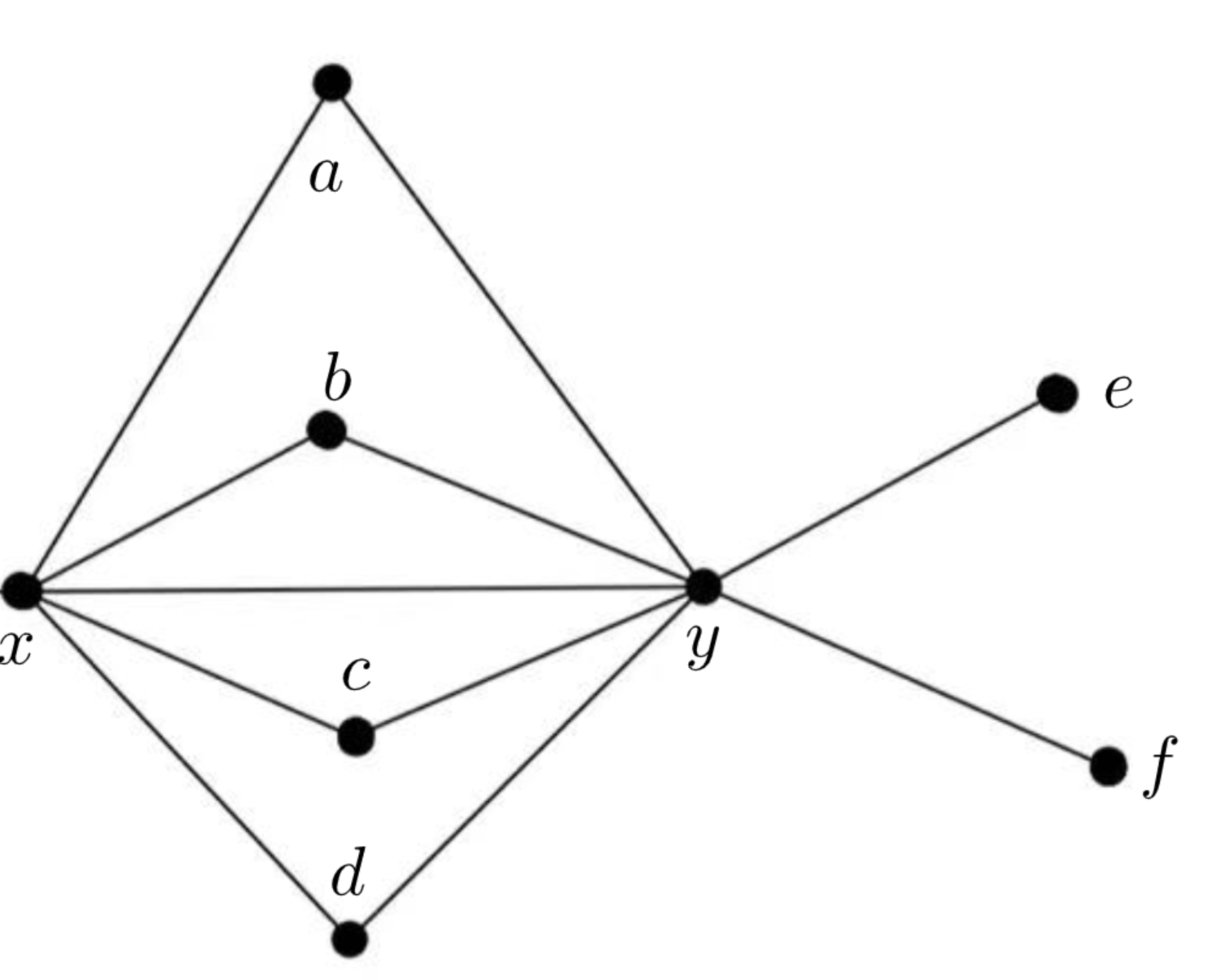}
    \caption{The 5-7 edge $xy$ is contained in 4 triangles.}
    \label{Fig:57354} 
\end{figure}

Removing the vertices in $H$ from the graph $G$. Assume that $a$ has two neighbors in $V(G)\backslash H$, we immediately get an $S_{3,5}$. Thus, any vertices in $S_1$ can have at most one neighbor in $V(G)\backslash H$.  
It can form a path of length 3 within $S_1$ and the edges of $S_1$ do not affect the neighbors of the vertices of $S_1$ in $V(G)\backslash H$. Each of the vertices $e$ and $f$ can have at most two neighbors in $V(G)\backslash H$. If a vertex $s$ in $S_2$ is adjacent to a vertex $t$ in $S_1$, this adjacent situation does not alter the set of vertices adjacent to $t$ within $S_1$ in the subgraph $V(G)\backslash H$. If $e$ and $f$ are adjacent, this implies that each of them has at most one neighbor in $V(G)\backslash H$. If $e$ is adjacent to at least one vertex in $S_1 \cup \{f\}$, then it has at most one neighbor in $V(G)\backslash H$. Similarly, the same situation holds for $f$.

Here, $y$ is a vertex of degree 7 in $G$. According to Claim 2.6, it is only necessary to consider the case where all the vertices in $S_3$ have at least one other neighbor within $S_3$. Specifically, consider vertices $e$ and $f$ in $S_3$.

Removing the vertices in $H$ from the graph $G$, the number of edges deleted is at most $11+3+4\times1+(1+2)\times2$=24. Similarly, a planar graph on eight vertices has at most $3\times8-6$=18 edges and the set $S_1\cup S_2$ contains vertices with at most one neighbor in $V(G)\backslash H$, there exist at most six additional edges beyond these 18 edges. At this case, both $e$ and $f$ are adjacent to at least two other vertices in $S_1\cup S_2$, respectively. Then the two vertices $e$ and $f$ are connected to at most six edges in total. If we delete 24 edges, then all the vertices in $S_1\cup S_2$ have at least one neighbor in $V(G)\backslash H$. At this case, there are two possible cases, as shown in the Figure 8.

\begin{figure}[h]
    \centering
    \includegraphics[width=0.6\textwidth]{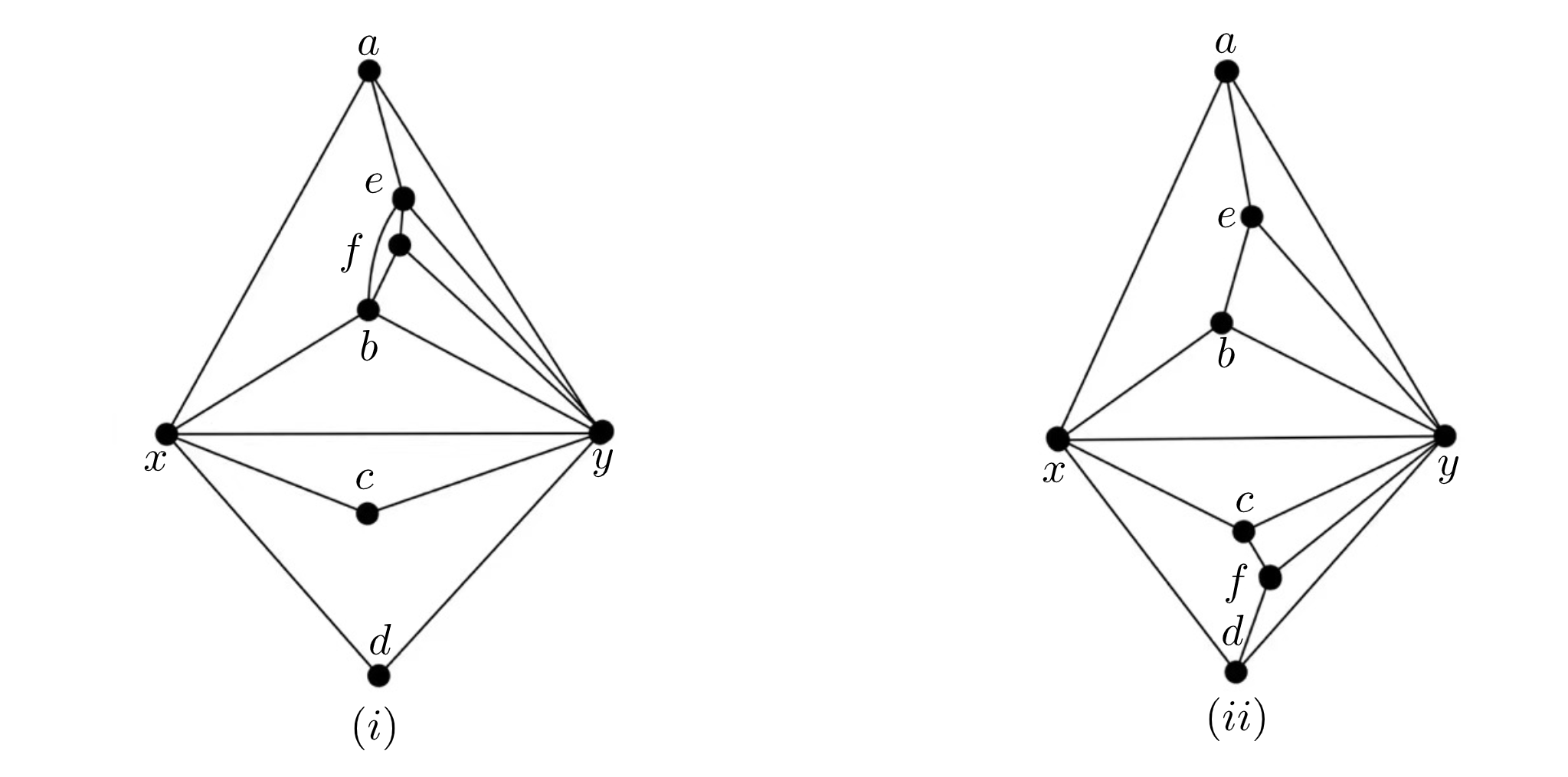}
    \caption{If we delete 24 edges, here are the adjacent situations of the vertices of $H$.}
    \label{Fig:57444} 
\end{figure}

For the Figure 8(i), at this case, the vertices $e$ and $f$ lie in the same face. When deleting 24 edges, there exists an edge cut set with size at most 2. Therefore, in this case, we must delete at most 23 edges to avoid finding an edge cut set with size at most 4.
 
For the Figure 8(ii), the vertices $e$ and $f$ lie in the different faces. If we delete 24 edges, there exists an edge cut set with size at most 3. Therefore, in this case, we must delete at most 23 edges to avoid finding an edge cut set with size at most 4.

Therefore, by the induction hypothesis,, Claim 2.9 holds. 

\textbf{Case 2. There are 3 triangles containing the edge $xy$.} 

Let $a$, $b$ and $c$ be the vertices in $G$ that are adjacent to both $x$ and $y$. Let $d$, $e$ and $f$ be the
vertices adjacent to $y$ but not to $x$, and $m$ be the
vertex adjacent to $x$ but not to $y$, as shown in the Figure 9. Define $S_1$=\{$a$, $b$, $c$\}, $S_2$=\{$d$, $e$, $f$\}, $S_3$=$S_1\cup S_2\cup\{x\}$ and let $H$=$S_3\cup\{y,m\}$. 

\begin{figure}[h]
    \centering
    \includegraphics[width=0.4\textwidth]{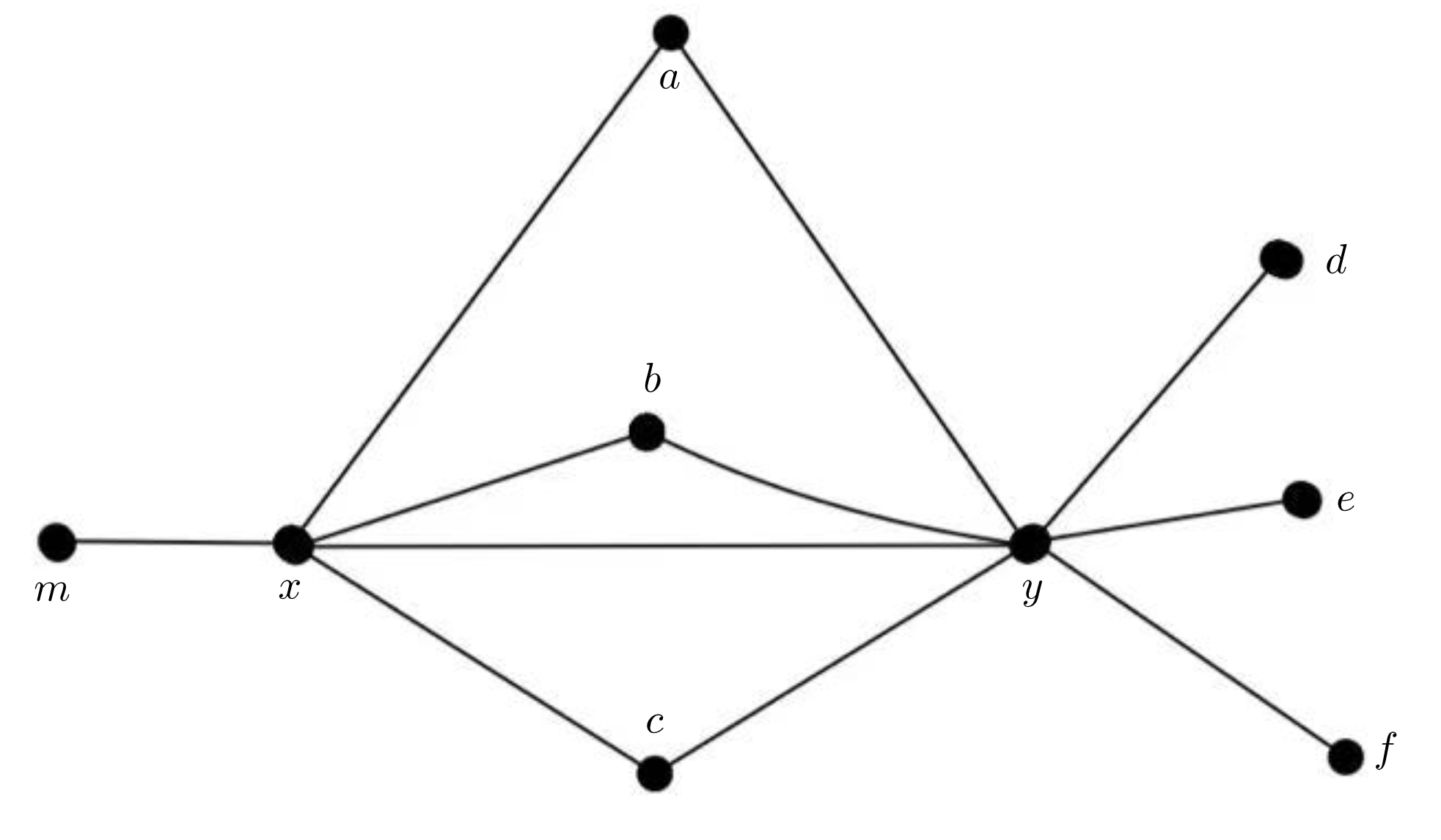}
    \caption{The 5-7 edge $xy$ is contained in 3 triangles.}
    \label{Fig:57353} 
\end{figure}

Removing the vertices in $H$ from the graph $G$. Assume that $a$ has two neighbors in $V(G)\backslash H$, we immediately get an $S_{3,5}$. Thus, any vertices in $S_1$ can have at most one neighbor in $V(G)\backslash H$.  
It can form a path of length 2 within $S_1$ and the edges of $S_1$ do not affect the neighbors of the vertices of $S_1$ in $V(G)\backslash H$. Each vertex in $S_2$ can have at most two neighbors in $V(G)\backslash H$.
If a vertex $r$ in $S_2$ is adjacent to a vertex $p$ in $S_1$, this adjacent situation does not alter the set of vertices adjacent to $p$ within $S_1$ in the subgraph $V(G)\backslash H$. If a vertex $r$ in $S_2$ is adjacent to at least one vertex in $S_1 \cup S_2 $, it makes these vertices in $S_2$ have at most one neighbor in $V(G)\backslash H$. The vertex $m$ can have at most four neighbors in $V(G)\backslash H$. If $m$ is adjacent to a vertex in $S_1$, then this vertex in $S_1$ have no neighbors in $V(G)\backslash H$, and $m$ has at most three neighbors in $V(G)\backslash H$. If $m$ is adjacent to a vertex in $S_2$, then this vertex in $S_2$ has at most one neighbor in $V(G)\backslash H$, and $m$ has at most three neighbors in $V(G)\backslash H$.

Here, $y$ is a vertex of degree 7. According to Claim 2.6, it is only necessary to consider the case where all the vertices in $S_3$ have at least one neighbor within $S_3$. Specifically, consider vertices $d$, $e$ and $f$ in $S_3$.

Deleting the vertices in $H$ from the graph $G$, the number of edges deleted is at most $11+2+3\times1+(1+2)\times3+4$=29. 

For any planar embedding of the vertex set $S_3\cup x$ (with $|S_3\cup x|=8$), the edge count is bounded by $3\times8-6$=18 due to Euler's formula. Given that each vertex in $S_1\cup S_2$ has at most one neighbor in $V(G)\backslash H$, there are at most six additional edges beyond these 18 edges. Furthermore, vertex $m$ has at most five incident edges. Therefore, the maximum number of edges removed is bounded by 29. At this case, $d$, $e$, and $f$ are adjacent to at least two other vertices in $S_1\cup S_2$, respectively. And the three vertices $d$, $e$, and $f$ are connected to at most nine edges in total. 

If the three vertices $d$, $e$, and $f$ lie in the same face and we deleting 29 edges, it is impossible that none of the three vertices $d$, $e$, and $f$ is adjacent to each other. The possible cases as shown in Figure 10.

\begin{figure}[h]
    \centering
    \includegraphics[width=0.9\textwidth]{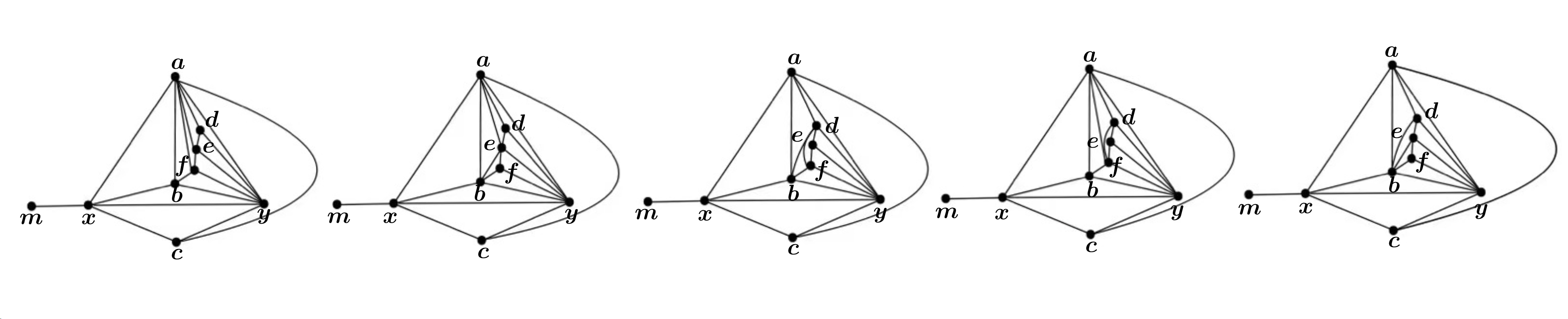}
    \caption{The vertices $d$, $e$, and $f$ lie in the same face.
}
    \label{Fig:5735} 
\end{figure}

If the three vertices $d$, $e$, and $f$ are not in the same faces, then there are two vertices among $d$, $e$, and $f$ that are adjacent to two vertices in $S_1$. Assume that these two vertices are $d$ and $e$. However, in a planar graph, the neighbors of $d$ and $e$ in $S_1$ cannot be exactly the same. Without loss of generality, let $d$ be adjacent to $a$ and $b$, and $e$ be adjacent to $c$ and $d$. At this case, $f$ is adjacent to one vertex in $S_1$ and one vertex in $\{d,e\}$, respectively. 
There are two cases as shown in the Figure 11. At this case, $m$ is not adjacent to any vertices in $S_1\cup S_2$.

\begin{figure}[h]
    \centering
    \includegraphics[width=0.6\textwidth]{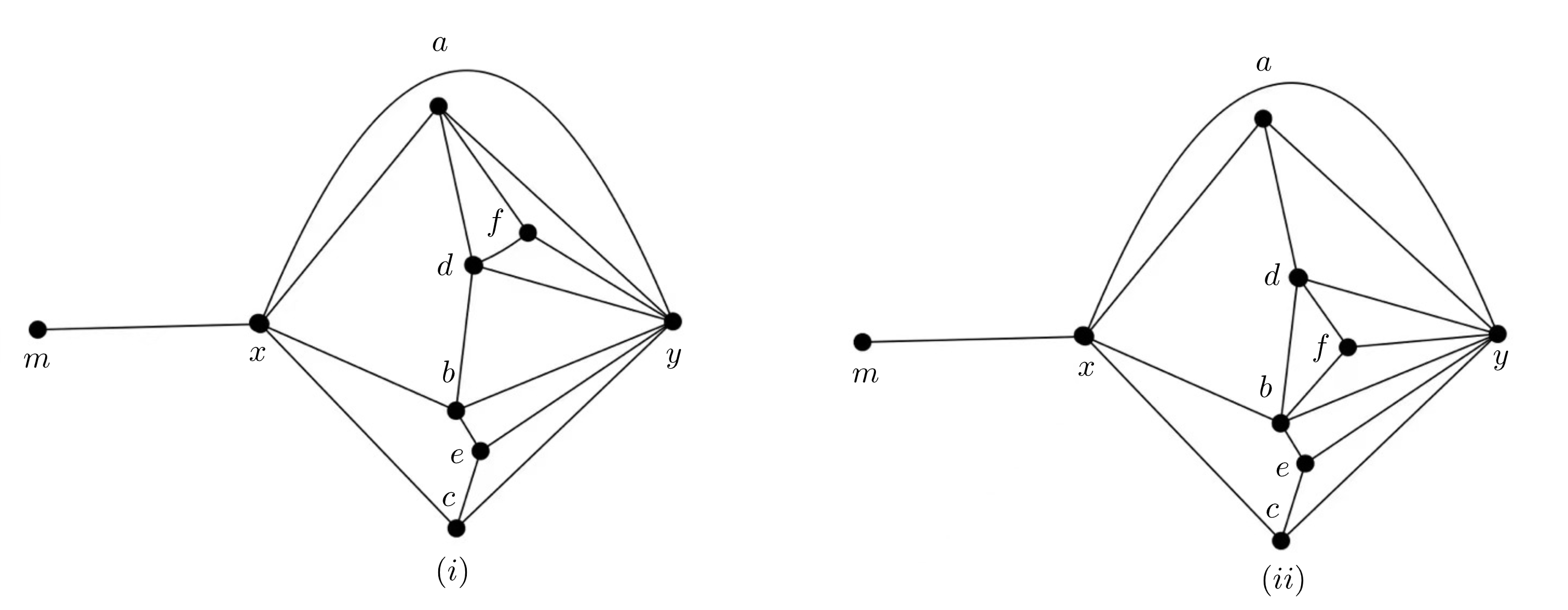}
    \caption{The vertices $d$, $e$, and $f$ lie in the different faces.
}
    \label{Fig:57222} 
\end{figure}

 For the Figure 10 or Figure 11, if we delete 29 edges, there must exist an edge cut set with size at most 4. So in any cases, we can delete at most 28 edges. 
 
When removing between 13 edges and 28 edges, there must exist a connected component $H_l$ in the resulting graph $G\backslash G[H] $, because there exist the vertex in $S_1\cup S_2$ have at least one neighbor in $V(G)\backslash H$.

Suppose that there are $l$ vertices in $H_l$ for $l\geq1$. If $l=1$, then we can  removing the vertices of $H\cup H_l$, there are 10 vertices and 28 edges. By the induction hypothesis, Claim 2.9 holds.    

Otherwise $l\geq2$, $H_l$ is a subgraph of the graph $G$. Therefore, $H_l$ does not contain any 6-6 edges, 5-6 edges or $S_{3,5}$ as a subgraph. Because the number of vertices in $H_l$ is less than $n$. By the induction hypothesis, $H_l$ has at most $\tfrac{23l}{8}-\tfrac{9}{2}$ edges. Let $V_l$ be the vertex set of $H_l$. Then, if we delete $H \cup V_l $ from the graph $G$, we will delete at most $\tfrac{23l}{8}-\tfrac{9}{2}$+28 edges.

By the induction hypothesis,

$e(G\backslash(H \cup V_l))$ $\leq$ $\tfrac{23(n-(l+9))}{8}-\tfrac{9}{2}$+$\tfrac{23l}{8}-\tfrac{9}{2}$+28
	$\leq$$\tfrac{23n}{8}-\tfrac{9}{2}$.

Therefore, deleting at most 28 edges, Claim 2.9 holds. 
$\hfill\square$ 

\vspace{3mm}
\noindent\textbf{Claim 2.10. } If the graph $G$ contains a 6-7 edge, then $e(G)$  $\leq$ $\tfrac{23n}{8}-\tfrac{9}{2}$.

\noindent\textbf{Proof.} Let $xy$$\in$$E(G)$ be a 6-7 edge. There are at least 4 triangles containing the edge $xy$. Otherwise, $G$ contains an $S_{3,5}$. We subdivide the cases based on the number of triangles containing the edge $xy$.

\textbf{Case 1. There are 5 triangles containing the edge $xy$.} 

Let $a$, $b$, $c$, $d$ and $e$ be the vertices in $G$ that are adjacent to both $x$ and $y$. Let $f$ be the
vertices adjacent to $y$ but not to $x$, as shown in the Figure 12. Define $S_1$=\{$a$, $b$, $c$, $d$, $e$\},$S_2$=$S_1\cup\{f,x\}$ and let $H$=$S_2\cup\{y\}$.

\begin{figure}[h]
    \centering
    \includegraphics[width=0.3\textwidth]{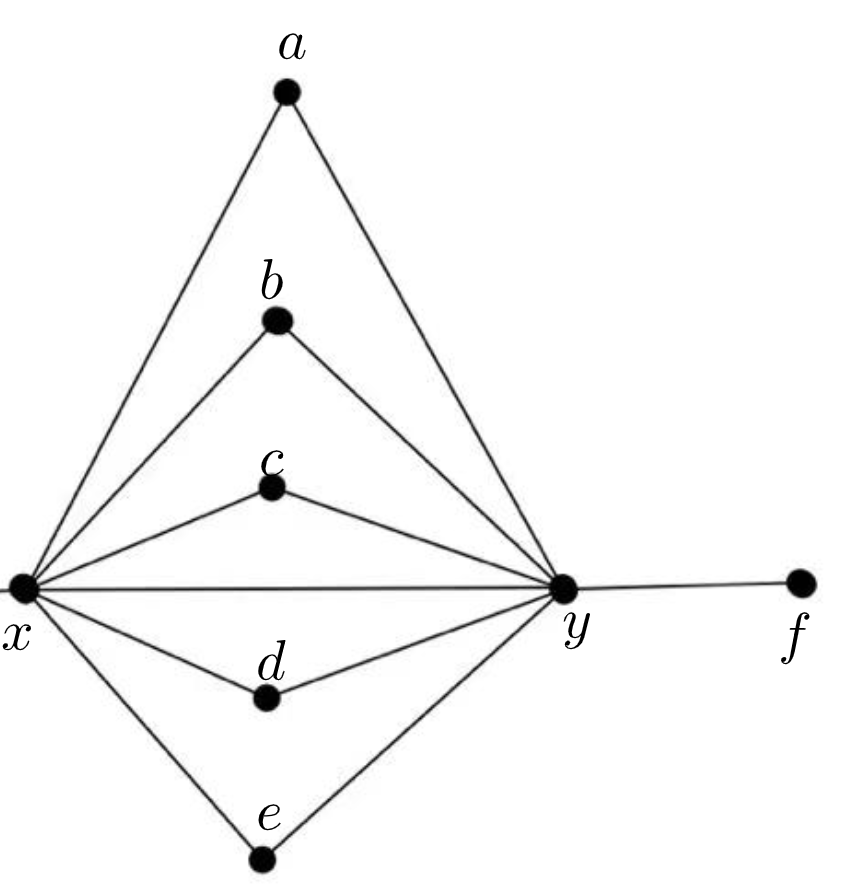}
    \caption{The 6-7 edge $xy$ is contained in 5 triangles.}
    \label{Fig:67355} 
\end{figure}

Deleting the vertices in $H$ from the graph $G$. Assume that $a$ has two neighbors in $V(G)\backslash H$, we immediately get an $S_{3,5}$. Thus, any vertices in $S_1$ can have at most one neighbor in $V(G)\backslash H$.  
It can form a path of length 4 within $S_1$ and the edges of $S_1$ do not affect the neighbors of the vertices of $S_1$ in $V(G)\backslash H$. The vertex $f$ can have at most two neighbors in $V(G)\backslash H$.
If $f$ is adjacent to a vertex $p$ in $S_1$, this adjacent situation does not alter the set of vertices adjacent to $p$ within $S_1$ in the subgraph $V(G)\backslash H$. If $f$ is adjacent to at least one vertex in $S_1$, then $f$ has at most one neighbor in $V(G)\backslash H$.

Here, $y$ is a vertex of degree 7. According to Claim 2.6, it is only necessary to consider the case where all the vertices in $S_2$ have at least one neighbor within $S_2$. Specifically, consider vertex $f$ in $S_2$.

Deleting the vertices in $H$ from the graph $G$, the number of edges deleted is at most $12+2+4+5\times1+1$=24. 

For any planar embedding of the vertex set $H$, the edge count is bounded by $3\times8-6$=18 due to Euler's formula. Given that each vertex in $S_1\cup {f}$ has at most one neighbor in $V(G)\backslash H$, there are at most six additional edges beyond these 18 edges. At this case, $f$ is adjacent to two vertices in $S_1$. If the number of deleted edges is 24, each vertex of $S_2\backslash\{x\}$ has one neighbor in $V(G)\backslash H$. In this case, there exists an edge cut set with size at most 2. Therefore, it is impossible to delete 24 edges. The number of edges deleted is at most 23.

By the induction hypothesis, 
$e(G\backslash H)$ $\leq$  $\tfrac{23(n-8)}{8}-\tfrac{9}{2}$,
    
    $e(G)$$\leq$ $e(G\backslash H)$+$23$ 
    $\leq$  $\tfrac{23(n-8)}{8}-\tfrac{9}{2}+23$ 
    $\leq$  $\tfrac{23n}{8}-\tfrac{9}{2}$.

Claim 2.10 is true.

\textbf{Case 2. There are 4 triangles containing the edge $xy$.} 

Let $a$, $b$, $c$ and $d$ be the vertices in $G$ that are adjacent to both $x$ and $y$. Let $e$ and $f$ be the
vertices adjacent to $y$ but not to $x$, $m$ be the
vertex adjacent to $y$ but not to $x$, as shown in Figure 13. Define $S_1$=\{$a$, $b$, $c$, $d$\}, $S_2$=\{$e$, $f$\}, $S_3$=$S_1\cup S_2\cup\{x\}$ and let $H$=$S_3\cup\{y,m\}$.

\begin{figure}[h]
    \centering
    \includegraphics[width=0.55\textwidth]{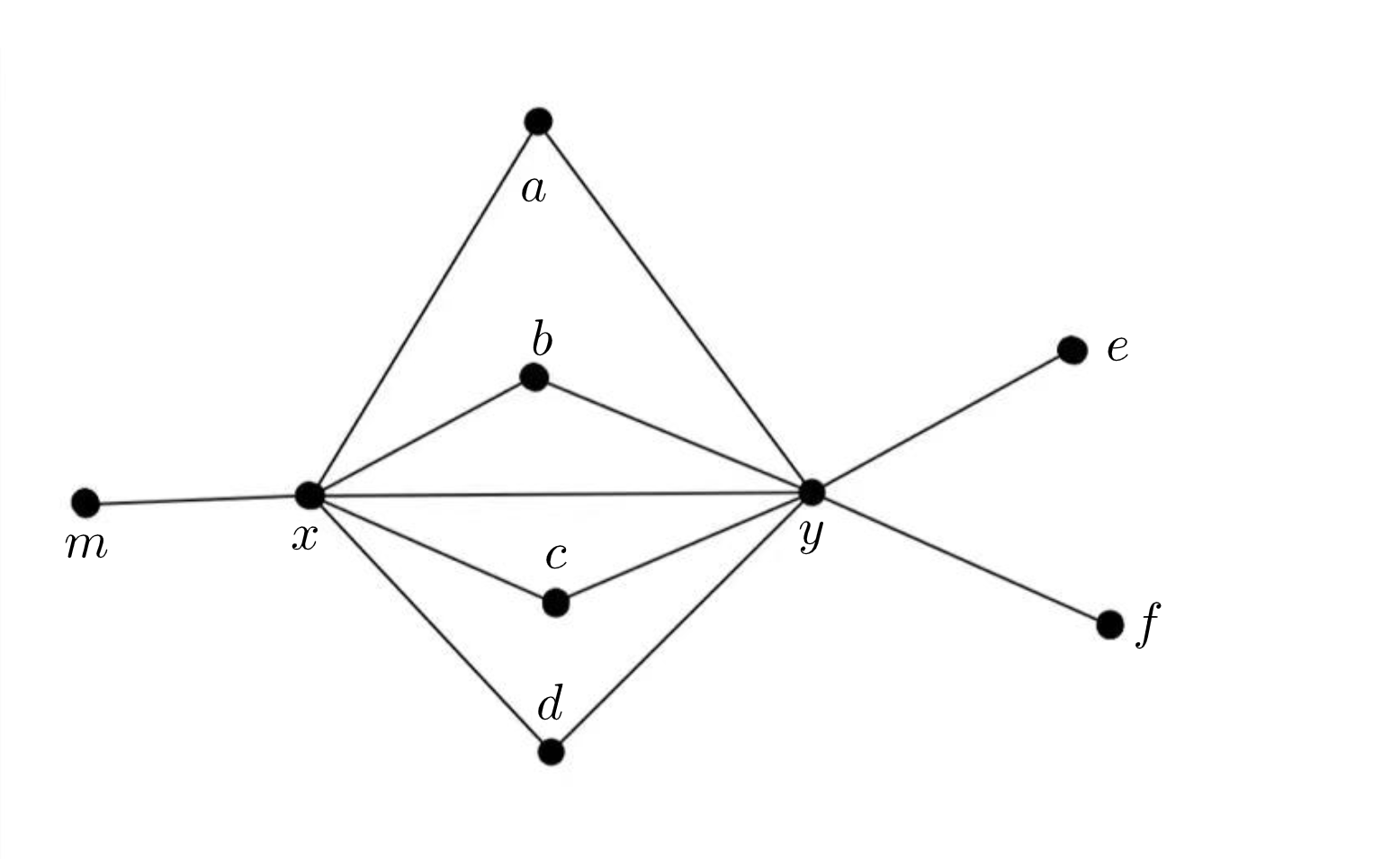}
    \caption{The 6-7 edge $xy$ is contained in 4 triangles.}
    \label{Fig:67354} 
\end{figure}

Deleting the vertices in $H$ from the graph $G$. Assume that $a$ has two neighbors in $V(G)\backslash H$, we immediately get an $S_{3,5}$. Thus, any vertices in $S_1$ can have at most one neighbor in $V(G)\backslash H$.  
It can form a path of length 3 within $S_1$ and the edges of $S_1$ do not affect the neighbors of the vertices of $S_1$ in $V(G)\backslash H$. Each vertex of $S_2$ can have at most two neighbors in $V(G)\backslash H$. If a vertex $q$ in $S_1$ is adjacent to a vertex $r$ in $S_2$, then the vertex $q$ in $S_1$ has at most one neighbor in $V(G)\backslash H$. If a vertex $t$ in $S_1$ is adjacent to two vertices $t_1$ and $t_2$ in $S_2$, then the vertex $t$ in $S_1$ has no neighbors in $V(G)\backslash H$. If a vertex $t_3$ in $S_2$ is adjacent to at least one vertex in $S_1 \cup S_2 $, it makes the vertex $t_3$ in $S_2$ has at most one neighbor in $V(G)\backslash H$. The vertex $m$ has at most two neighbors in $V(G)\backslash H$. If $m$ is adjacent to a vertex $t_4$ in $S_1$, then this vertex $t_4$ in $S_1$ has no neighbors in $V(G)\backslash H$, and $m$ has at most two neighbors in $V(G)\backslash H$. If $m$ is adjacent to a vertex $t_5$ in $S_2$, then the vertex $t_5$ in $S_2$ has at most one neighbor in $V(G)\backslash H$, and $m$ has at most one neighbor in $V(G)\backslash H$ as well.

Here, $y$ is a vertex of degree 7. According to Claim 2.6, it is only necessary to consider the case where all the vertices in $S_3$ have at least one neighbor within $S_3$. Specifically, consider the vertices $e$ and $f$ in $S_3$.

Deleting the vertices in $H$ from the graph $G$, the number of edges deleted is at most $12+3+2+4\times1+(1+2)\times2$=27. 

For any planar embedding of the vertex set $S_3\cup {y}$ (with $|S_3\cup {y}|=8$), the edge count is bounded by $3\times8-6$=18 due to Euler's formula. Given that each vertex in $S_1\cup S_2$ has at most one neighbor in $V(G)\backslash H$, there are at most six additional edges beyond these 18 edges. Furthermore, vertex $m$ has at most three incident edges. Therefore, the maximum number of edges removed is bounded by 27.
At this case,
both $e$ and $f$ are adjacent to at least two other vertices in $S_1\cup S_2$, respectively.
In any cases, apart from the two edges $ey$ and $ef$, the two vertices $e$ and $f$ are connected to four edges in total.

At this case, the adjacent situation of the vertices in $H$ is as shown in the Figure 14.

\begin{figure}[h]
    \centering
    \includegraphics[width=0.7\textwidth]{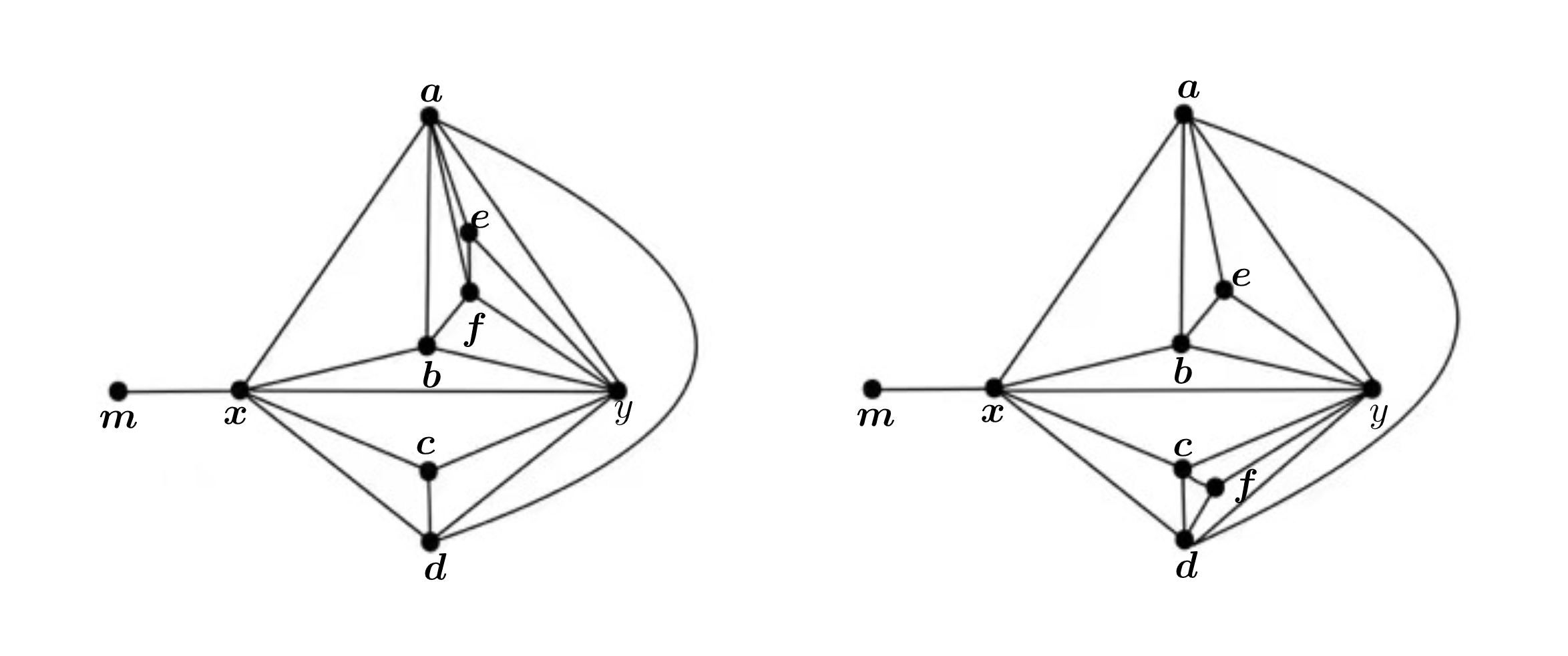}
    \caption{If we delete 27 edges, here are the adjacent situations of the vertices of $H$.}
    \label{Fig:67333} 
\end{figure}

 when removing between 15 edges and 27 edges, there must exist a connected component $H_l$ in the resulting graph $G\backslash G[H] $, because there exist the vertex in $S_1\cup S_2$ have at least one neighbor in $V(G)\backslash H$.

Suppose that there are $l$ vertices in $H_l$ for $l\geq1$. If $l=1$, then we can  removing the vertices of $H\cup H_l$, there are 10 vertices and 27 edges. By the induction hypothesis, Claim 2.10 holds.    

Otherwise $l\geq2$, $H_l$ is a subgraph of the graph $G$. Therefore, $H_l$ does not contain any 6-6 edges, 5-6 edges or $S_{3,5}$ as a subgraph. Because the number of vertices in $H_l$ is less than $n$. By the induction hypothesis, $H_l$ has at most $\tfrac{23l}{8}-\tfrac{9}{2}$ edges. Let $V_l$ be the vertex set of $H_l$. Then, if we delete $H \cup V_l $ from the graph $G$, we will delete at most $\tfrac{23l}{8}-\tfrac{9}{2}$+27 edges.

By the induction hypothesis,

$e(G\backslash(H \cup V_l))$ $\leq$ $\tfrac{23(n-(l+9))}{8}-\tfrac{9}{2}$+$\tfrac{23l}{8}-\tfrac{9}{2}$+27
	$\leq$$\tfrac{23n}{8}-\tfrac{9}{2}$.

Therefore, deleting at most 27 edges, Claim 2.10 holds. 
$\hfill\square$ 

\vspace{3mm}
\noindent\textbf{Claim 2.11. } If the graph $G$ contains a 7-7 edge, then $e(G)$  $\leq$ $\tfrac{23n}{8}-\tfrac{9}{2}$.

\noindent\textbf{Proof.} Let $xy$$\in$$E(G)$ be a 7-7 edge. There are at least 5 triangles containing the edge $xy$. Otherwise, $G$ contains an $S_{3,5}$. We subdivide the cases based on the number of triangles containing the edge $xy$.

\textbf{Case 1. There are 6 triangles containing the edge $xy$.} 

Let $a$, $b$, $c$, $d$, $e$ and $f$ be the vertices in $G$ that are adjacent to both $x$ and $y$, as shown in the Figure 15. Define $S_1$=\{$a$, $b$, $c$, $d$, $e$, $f$\}, $S_2$=$S_1\cup\{x\}$ and $H$=$S_2\cup\{y\}$. 

\begin{figure}[h]
    \centering
    \includegraphics[width=0.3\textwidth]{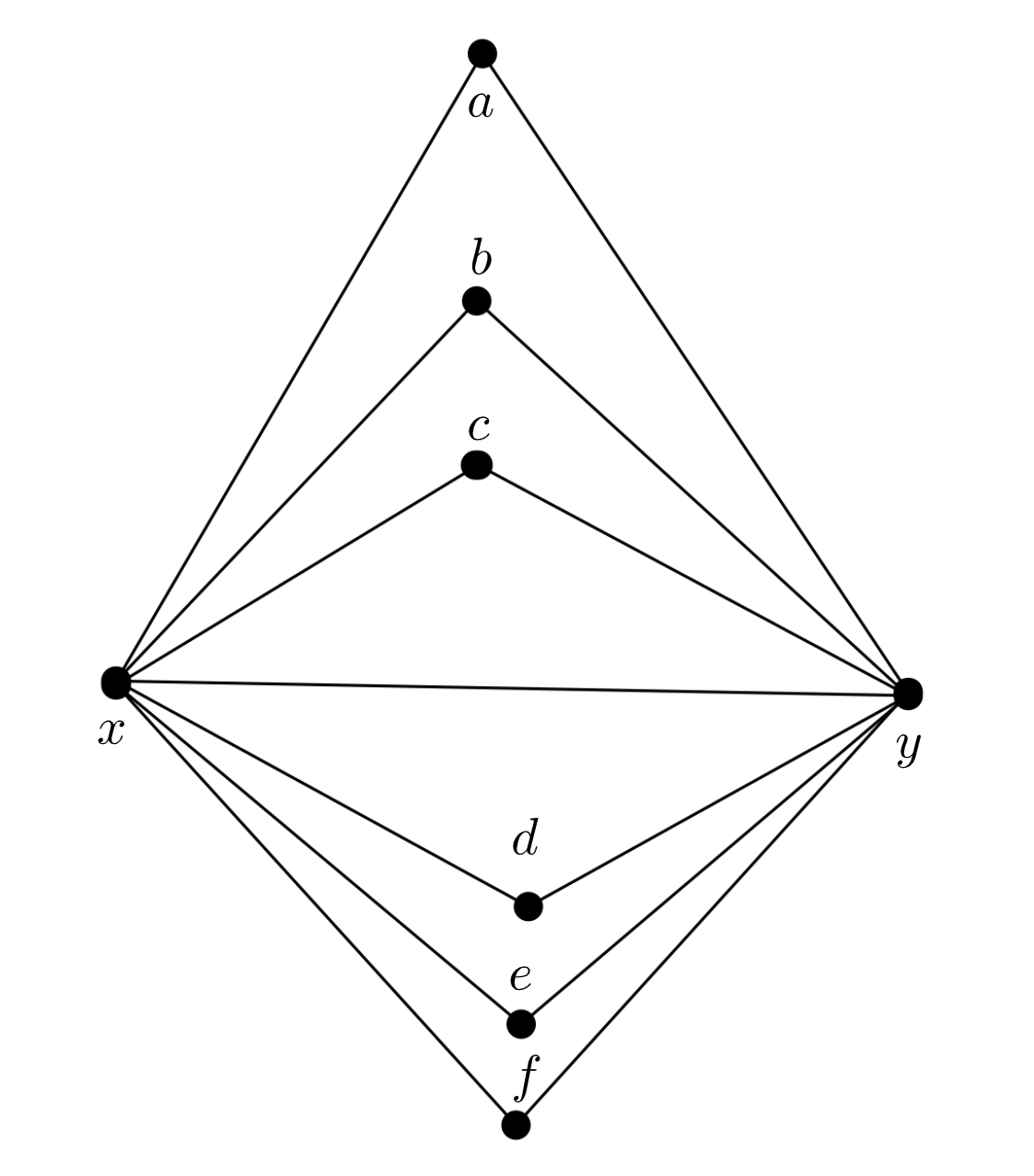}
    \caption{The 7-7 edge $xy$ is contained in 6 triangles.}
    \label{Fig:77356} 
\end{figure}

Deleting the vertices in $H$ from the graph $G$. Assume that $a$ has two neighbors in $V(G)\backslash H$, then we immediately get an $S_{3,5}$. Thus, any vertices in $S_1$ can have at most one neighbor in $V(G)\backslash H$.  
It can form a path of length 5 within $S_1$ and the edges of $S_1$ do not affect the neighbors of the vertices of $S_1$ in $V(G)\backslash H$.
The number of edges deleted is at most $13+5+6\times1+$=24.

For any planar embedding of the vertex set $H$, the edge count is bounded by $3\times8-6$=18 due to Euler's formula. Given that each vertex in $S_1$ has at most one neighbor in $V(G)\backslash H$, there are at most six additional edges beyond these 18 edges. If we delete 24 edges, each vertex in $S_1$ has a neighbor in $V(G)\backslash H$. At this case, there exists an edge cut set with size at most 2, a contradiction.

Deleting the vertices in $H$ from the graph $G$, the number of edges deleted is at most 23. By the induction hypothesis, the Claim 2.11 is true.

\textbf{Case 2. There are 5 triangles containing the edge $xy$.} 

Let $a$, $b$, $c$, $d$ and $e$ be the vertices in $G$ that are adjacent to both $x$ and $y$. Let $g$ be the
vertices adjacent to $y$ but not to $x$, $f$ be the
vertex adjacent to $x$ but not to $y$, as shown in the Figure 16. Define $S_1$=\{$a$, $b$, $c$, $d$, $e$\}, $S_2$=$S_1\cup \{g,x\}$ and let $H$=$S_1\cup \{f,y\}$. 

\begin{figure}[h]
    \centering
    \includegraphics[width=0.35\textwidth]{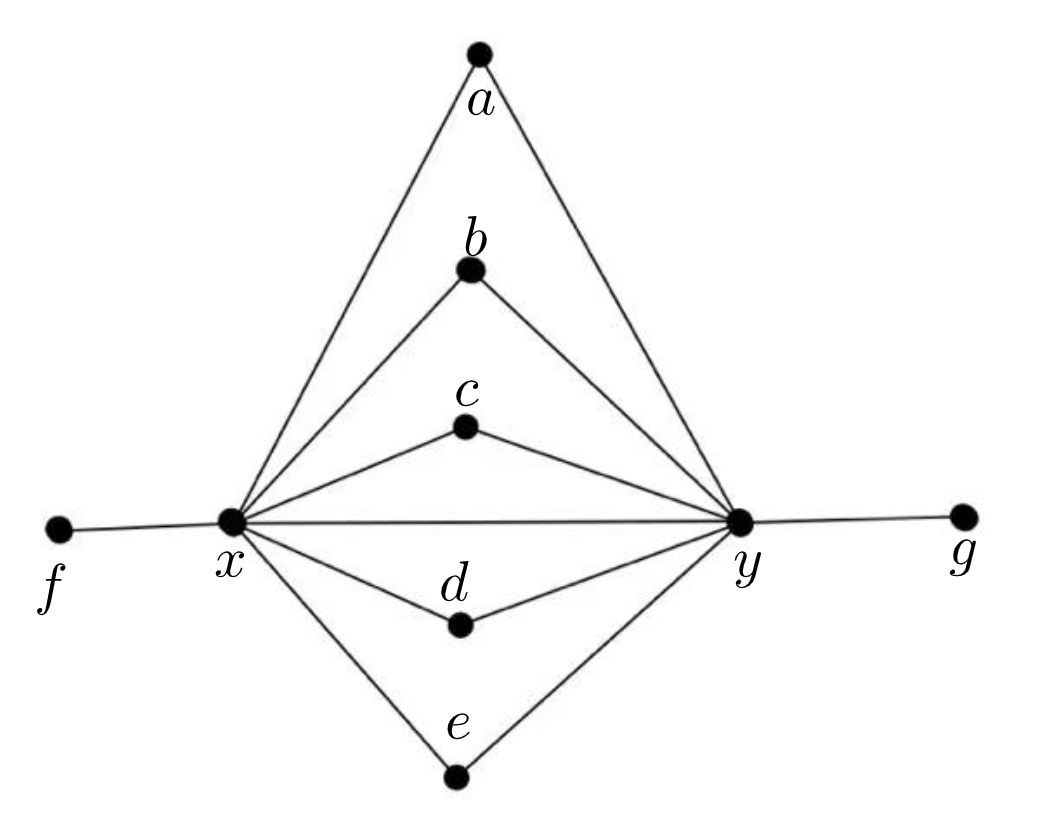}
    \caption{The 7-7 edge $xy$ is contained in 5 triangles}
    \label{Fig:77355} 
\end{figure}

Deleting the vertices in $H$ from the graph $G$. Assume that $a$ has two neighbors in $V(G)\backslash H$, then we immediately get an $S_{3,5}$. Thus, any vertices in $S_1$ can have at most one neighbor in $V(G)\backslash H$.  
It can form a path of length 4 within $S_1$ and the edges of $S_1$ do not affect the neighbors of the vertices of $S_1$ in $V(G)\backslash H$. Both $f$ and $g$ can have at most two neighbors in $V(G)\backslash H$. If $f$ (or $g$) is adjacent to at least one vertex in $S_1$, then these vertices in $S_1$ have no neighbors in $V(G)\backslash H$, and $f$ (or $g$) has at most one neighbor in $V(G)\backslash H$. If $f$ is adjacent to $g$, then the maximum number of neighbors of $f$ and $g$ in $V(G)\backslash H$ will decrease by 1 each.

Here, $y$ is a vertex of degree 7. According to Claim 2.6, it is only necessary to consider the case where all the vertices in $S_2$ have at least one neighbor within $S_2$. Specifically, consider the vertex $g$ in $S_2$.

Deleting the vertices in $H$ from the graph $G$, the number of edges deleted is at most $13+4+4+2\times1+3\times1$=26.

At this case, the two neighbors of $f$ and $g$ in $S_1$ are the same. If we delete 26 edges, there exists an edge cut set with size at most 2, a contradiction.

Deleting the vertices in $H$ from the graph $G$, the number of edges deleted is at most 25. By the induction hypothesis, the Claim 2.11 is true.$\hfill\square$

\section{Proof of Theorem 1.1.}
\noindent\textbf{Theorem 1.1.}  Let $G$ be an $S_{3,5}$-free planar graph with $n$ vertices. Then $ex_\mathcal{P}(n, S_{3,5})$ $\leq$ $\tfrac{23n}{8}-\tfrac{9}{2}$ for all $n \geq 2$.

\noindent\textbf{Proof.} Let $xy$ be an edge of the graph $G$. If $d(x)+d(y)$$\geq$15, then $\Delta(G)$ $\geq 8$. Theorem 1.1 holds by Claim 2.4 and Claim 2.5. 

 If $d(x)+d(y)$=14, and $\Delta(G)$ $\leq 8$, then $xy$ must be a 7-7 edge. Theorem 1.1 holds by Claim 2.11.

 If $d(x)+d(y)$=13, and $\Delta(G)$ $\leq 8$, then $xy$ must be a 6-7 edge. Theorem 1.1 holds by Claim 2.10.

 If $d(x)+d(y)$=12, and $\Delta(G)$ $\leq 8$, then $xy$ must be a 5-7 edge or a 6-6 edge. Theorem 1.1 holds by Claim 2.9 and Claim 2.7.

If $d(x)+d(y)$=11, and $\Delta(G)$ $\leq 8$, then $xy$ must be a 4-7 edge or a 5-6 edge. Theorem 1.1 holds by Claim 2.8 and Claim 2.7.
 
If $d(x)+d(y)$  $\leq $10, then 
 
\noindent $\sum_{xy \in E(G)}$ $(d(x)+d(y))$=$\sum_{xy \in E(G)}$ $d(x)^2$ $\geq$ $n\cdot \bar{d}^2$=$n\cdot  (\frac{2e}{n})^2 $, where $e\leq \frac{5n}{2} \leq \tfrac{23n}{8}-\tfrac{9}{2}$ for any $n \geq 12$.

 Combined with Claim 2.1  we know that Theorem 1.1 holds.$\hfill\square$

\section*{Funding}
This work was funded in part by Natural Science Foundation of China (Grant No. 12071194).

\section*{Data availability}
No data was used for the research described in the article.

\bibliographystyle{plain}
\bibliography{main}

\begin{thebibliography}{10}

\bibitem{Dowden2015ExtremalCP}
Chris Dowden.
\newblock {Extremal $C_4$‐free/$C_5$‐free planar graphs}.
\newblock {\em Journal of Graph Theory}, 83, 2015.

\bibitem{Du2021PlanarTN}
Lian~Yan Du, Bing Wang, and Mingqing Zhai.
\newblock {Planar Turán numbers on short cycles of consecutive lengths}.
\newblock {\em Bulletin of the Iranian Mathematical Society}, 48:2395 -- 2405, 2021.

\bibitem{Fang2023ExtremalSR}
Longfei Fang, Huiqiu Lin, and Yongtang Shi.
\newblock Extremal spectral results of planar graphs without vertex‐disjoint cycles.
\newblock {\em Journal of Graph Theory}, 106:496 -- 524, 2023.

\bibitem{2021arXiv211010515G}
Debarun {Ghosh}, Ervin {Gy{\H{o}}ri}, Addisu {Paulos}, and Chuanqi {Xiao}.
\newblock {Planar Turán number of double stars}.
\newblock arXiv:2110.10515, October 2021.

\bibitem{doi:10.1137/21M140657X}
Debarun Ghosh, Ervin Gy{\H{o}}ri, Ryan~R. Martin, Addisu Paulos, and Chuanqi Xiao.
\newblock {Planar Turán number of the 6-Cycle}.
\newblock {\em SIAM Journal on Discrete Mathematics}, 36(3):2028--2050, 2022.

\bibitem{2023planarturannumbersevencycle}
Ervin {Gy{\H{o}}ri}, Alan {Li}, and Runtian {Zhou}.
\newblock {The planar Turán number of the seven-cycle}.
\newblock arXiv:\\2307.06909, July 2023.

\bibitem{Lan2019PlanarTN}
Yongxin Lan and Yongtang Shi.
\newblock {Planar Turán numbers of short paths}.
\newblock {\em Graphs and Combinatorics}, 35:1035 -- 1049, 2019.

\bibitem{lan2018extremalhfreeplanargraphs}
Yongxin {Lan}, Yongtang {Shi}, and Zi-Xia {Song}.
\newblock {Extremal $H$-free planar graphs}.
\newblock arXiv:1808.\\01487, August 2018.

\bibitem{lan2019extremalthetafreeplanargraphs}
Yongxin Lan, Yongtang Shi, and Zi-Xia Song.
\newblock {Extremal Theta-free planar graphs}.
\newblock {\em Discrete Mathematics}, 342(12):111610, 2019.

\bibitem{shi2023planarturannumber7cycle}
Ruilin Shi, Zachary Walsh, and Xingxing Yu.
\newblock {Planar Turán number of the 7-cycle}.
\newblock {\em European Journal of Combinatorics}, 2023.

\bibitem{XU2024326}
Xin Xu, Yue Hu, and Xu~Zhang.
\newblock {An improved upper bound for planar Turán number of double star $S_{2,5}$}.
\newblock {\em Discrete Applied Mathematics}, 358:326--332, 2024.

\bibitem{XU2025114571}
Xin Xu and Jiawei Shao.
\newblock {The planar Turán number of double star $S_{2,4}$}.
\newblock {\em Discrete Mathematics}, 348(11):114571, 2025.

\bibitem{XU11}
Xin Xu, Xu~Zhang, and Jiawei Shao.
\newblock {Planar Turán number of double star $S_{3, 4}$}.
\newblock {\em AIMS Mathematics}, 10(1):1628--1644, 2025.

\end{thebibliography}
\end{document}